\documentclass[12pt,reqno]{amsart}
\usepackage{amsmath}
\usepackage{amssymb}
\usepackage{verbatim}
\usepackage{mathrsfs}
\usepackage{graphicx} 
\usepackage{caption}
\usepackage{subcaption}
\usepackage[noadjust]{cite}
\usepackage{hyperref}
\usepackage{color}

\usepackage[noadjust]{cite}

\usepackage[top=0.5in,bottom=0.6in,left=0.67in,right=0.67in]{geometry}

\newtheorem{lemma}{Lemma}
\newtheorem{prop}[lemma]{Proposition}

\newtheorem{theorem}[lemma]{Theorem}
\newtheorem{remark}[lemma]{Remark}
\newtheorem{rem}[lemma]{Remark}
\newtheorem{exmp}{Example}

\newtheorem{definition}[lemma]{Definition}

\numberwithin{equation}{section}
\numberwithin{lemma}{section}

\newcommand{\C}{\mathbb{C}}    
\newcommand{\N}{\mathbb{N}}    
\newcommand{\R}{\mathbb{R}}    
\newcommand{\Z}{\mathbb{Z}}    

\newcommand{\wh}{\widehat}
\renewcommand{\le}{\leqslant}
\renewcommand{\ge}{\geqslant}

\newcommand{\bo}{\mathscr{O}} 

\newcommand{\veps}{\varepsilon}

\newcommand{\DG}{{\mathsf{Diag}}}

\newcommand{\er}{\eqref}

\newcommand{\mra}{{\mathring{a}}}

\newcommand{\mrw}{{\mathring{w}}}

\newcommand{\mrphi}{\mathring{\phi}}

\newcommand{\pe}{\mathsf{e}}

\newcommand{\bp}{ \begin{proof} }
	\newcommand{\ep}{\hfill \end{proof} }
\newcommand{\be}{ \begin{equation} }
\newcommand{\ee}{ \end{equation} }
\newcommand{\tp}{\mathsf{T}}
\newcommand{\dm}{\mathsf{M}} 
\newcommand{\vgu}{\upsilon} 

\newcommand{\sr}{\operatorname{sr}}  
\newcommand{\sm}{\operatorname{sm}}  

\newcommand{\td}{\boldsymbol{\delta}}  

\newcommand{\sd}{\mathcal{S}}  

\newcommand{\cG}{\mathcal{G}}
\newcommand{\cH}{\mathcal{H}}

\newcommand{\ttm}{\tilde{m}}

\newcommand{\dR}{\mathbb{R}^d}

\newcommand{\dZ}{\mathbb{Z}^d}

\newcommand{\dLp}[1]{L_{#1}(\mathbb{R}^d)}

\newcommand{\dLrs}[3]{(L_{#1}(\mathbb{R}^d))^{#2\times #3}}

\newcommand{\fsupp}{\text{fsupp}}

\newcommand{\vertiii}[1]{{\left\vert\kern-0.25ex\left\vert\kern-0.25ex\left\vert #1
		 \right\vert\kern-0.25ex\right\vert\kern-0.25ex\right\vert}}

\newcommand{\dlp}[1]{l_{#1}(\mathbb{Z}^d)}

\newcommand{\dlpp}[1]{l_{#1}(\mathbb{Z}^2)}

\newcommand{\dlrs}[3]{(l_{#1}(\mathbb{Z}^d))^{#2\times #3}}

\newcommand{\vsp}[1]{\vspace{#1 cm}}

\newcommand{\drv}{\operatorname{drv}}
\newcommand{\mom}{\operatorname{mom}}
\newcommand{\spa}{\operatorname{span}}

\begin{document}

\title[Multivariate Vector Subdivisions] {Multivariate Vector Subdivision Schemes with a General Matrix-valued Filter}

\author{Ran Lu}
\address{School of Mathematics, Hohai University, Nanjing, China 211100.
	\quad  {\tt rlu3@hhu.edu.cn}}

\thanks{The research of the first author was supported by the National Natural Science Foundation of China under grant 12201178. }

	\makeatletter \@addtoreset{equation}{section} \makeatother
	
\begin{abstract}Subdivision schemes are closely related to splines and wavelets and have numerous applications in CAGD and numerical differential equations. Subdivision schemes employ a scalar filter; that is, scalar subdivision schemes, have been extensively studied in the literature. In contrast, subdivision schemes with a matrix filter, which are the so-called vector subdivision schemes, are far from being well understood. So far, only vector subdivision schemes that use special matrix-valued filters have been well-investigated, such as the Lagrange and Hermite subdivision schemes. To the best of our knowledge, it remains unclear how to define and characterize the convergence of a vector subdivision scheme that uses a general matrix-valued filter. Though filters from Lagrange and Hermite subdivision schemes have nice properties and are widely used in practice, filters not from either subdivision scheme appear in many applications. Hence, it is necessary to study vector subdivision schemes with a general matrix-valued filter. In this paper, from the perspective of a vector cascade algorithm, we show that there is only one meaningful way to define a vector subdivision scheme. We will analyze the convergence of the newly defined vector subdivision scheme and show that it is equivalent to the convergence of the corresponding vector cascade algorithm. Applying our theory, we show that existing results on the convergence of Lagrange and Hermite subdivision schemes can be easily obtained and improved. Finally, we will present some examples of vector subdivision schemes to illustrate our main results.

\end{abstract}


\maketitle

\pagenumbering{arabic}

\section{Introduction}

Subdivision schemes are fast-averaging algorithms to numerically compute refinable functions and their derivatives. Because of their multi-scale structure, subdivision schemes are closely related to splines and wavelets and have been numerously applied to computer-aided graphic design to generate smooth curves and surfaces (\cite{dd89,dl02,dlg87,dlg90,fhs22,hj98,hj06,hyx05,ms19}), numerically solving differential equations (\cite{hm21,hm23pp} and references therein), processing data with discrete wavelet and framelet transforms (\cite{hanbook,hl20pp}), and many other applications.

Most subdivision schemes that have been covered in the literature employ a scalar filter, which are called scalar subdivision schemes, and their properties and applications have been well-investigated (\cite{cdm91,cc13,dlg87,dlg90,dl02,hj06} and many references therein). By contrast, subdivision schemes that employ a matrix-valued filter, called vector (or matrix-valued) subdivision schemes, are much more complicated and are far from well-understood. So far, only vector subdivision schemes that use special matrix-valued filters have been well-studied, such as Lagrange subdivision schemes (named ``vector subdivision schemes" in \cite{dm06,hkz09,ms98}) and Hermite subdivision schemes and their generalizations \cite{ccmm21,ch19,cmss19,dm09,han22,han23-2,m92,ms17,ms19}. Though filters from Lagrange and Hermite subdivision schemes have nice properties and are widely used in practice, filters not from either scheme appear in many applications. For instance, matrix-valued filters that are of neither Lagrange nor Hermite type are used in \cite{cj05,cj06,cj09,js16} for generating curves and surfaces with specific properties;  also, many multiwavelets used in numerical PDEs are derived from refinable vector functions whose filters are of neither Lagrange or Hermite type (e.g.\cite{hm21,hm23pp}). Hence, it is necessary to work with vector subdivision schemes with an arbitrary matrix-valued filter and study their properties, such as convergence and smoothness. Unfortunately, to the best of our knowledge, it remains unclear how to define and characterize the convergence of a vector subdivision scheme that uses a general matrix-valued filter. 

In this paper, we introduce the notion of a vector subdivision scheme that works with a general matrix-valued filter and study its convergence. Since subdivision schemes are naturally related to cascade algorithms, we will perform our analysis of vector subdivision schemes from the perspective of vector cascade algorithms. Our study shows that there is only one meaningful way to define a vector subdivision scheme in that its convergence is equivalent to the convergence of the corresponding vector cascade algorithm.

\subsection{Cascade Algorithms and Subdivision Schemes}

We first recall the definition of a cascade algorithm. Let $r,s\in\N$, by an \emph{$r\times s$ matrix-valued $d$-dimensional filter} we mean a sequence $u=\{u(k)\}_{k\in\dZ}:\dZ\to \C^{r\times s}$ such that $u(k)\ne 0$ for only finitely many terms. By $\dlrs{0}{r}{s}$ we denote the linear space of all $r\times s$ matrix-valued $d$-dimensional filters. Let $\dm\in\N\setminus\{1\}$ be a dilation factor, $r\in\N$, $a\in\dlrs{0}{r}{r}$ and $p\in[1,\infty]$. Denote $I_d$ the $d\times d$ identity matrix, a \emph{vector $\dm I_d$-cascade operator with the filter $a$} on the space $\dLrs{p}{r}{1}$ is defined by
$$R_{a,\dm I_d}\psi(x)=\dm^d\sum_{k\in\dZ}a(k)\psi(\dm x-k),\quad\forall \psi\in(\dLp{p})^{r\times 1},\quad x\in\dR.$$
If $\phi$ is a fixed point of the operator $R_{a,\dm I_d}$, that is, 
\be\label{ref:eq}\phi(x)=R_{a,\dm I_d}\phi(x)=\dm^d\sum_{k\in\dZ}a(k)\phi(\dm x-k),\quad\forall x\in\dR,\ee 
then $\phi$ is called \emph{a vector $\dm I_d$-refinable vector function of the filter $a$} and in this case, $a$ is called an \emph{$\dm I_d$-refinement filter of $\phi$}. The equation \er{ref:eq} is called \emph{an $\dm I_d$-refinement equation} in the theory of subdivision schemes and wavelet analysis. When studying subdivision schemes, it is natural to write the refinement equation \er{ref:eq} in an equivalent form using the Fourier transform. For a matrix-valued function $f\in\dLrs{1}{s}{t}$, its \emph{Fourier transform} is given by 
$$\wh{f}(\xi):=\int_{\dR}f(x)e^{-ix\cdot\xi}dx,\quad\forall\xi\in\dR,$$
where the above integral is taken entry-wise. The definition of Fourier transform is naturally extended to matrices of $\dLp{2}$ functions or tempered distributions. For a filter $u\in\dlrs{0}{s}{t}$, define its \emph{Fourier series} by
$$\wh{u}(\xi):=\sum_{k\in\dZ}u(k)e^{-ik\cdot\xi},\quad\forall\xi\in\dR.$$
By taking the Fourier transform, the equation \er{ref:eq} is then equivalent to
\be\label{ref:eq:f}\wh{\phi}(\dm \xi)=\wh{a}(\xi)\wh{\phi}(\xi),\quad\forall\xi\in\dR.\ee
The algorithm that iteratively applies the operator $R_{a,\dm I_d}$ to $\psi\in\dLrs{p}{r}{1}$ and generates a sequence $\{R_{a,\dm I_d}^n\psi\}_{n=1}^\infty$ is called a \emph{vector $\dm I_d$-cascade algorithm} with the filter $a$. A cascade algorithm is often used to approximate or compute a refinable vector function and its derivatives. Suppose $\phi$ is an $\dm I_d$-refinable vector function of the filter $a$. If for some initial vector function $\psi$, the sequence $\{R_{a,\dm I_d}^n\psi\}_{n=1}^\infty$ converges to $\phi$ in some function space, then the cascade algorithm is said to be \emph{convergent}. The convergence of a cascade algorithm in $L_p$ or Sobolev spaces has been well-studied in \cite{han03,cjd02,li01,li03} and many references therein.

In scientific computing, we need a discretization of a cascade algorithm to compute a refinable vector $\phi$ and its derivatives. This motivates people to introduce a discrete analog to a cascade algorithm, which is a \emph{subdivision scheme}. Let $\dm\in\N\setminus\{1\}$, $r\in\N$ and $a\in\dlrs{0}{r}{r}$. A \emph{vector $\dm I_d$-subdivision operator with the filter $a$} is given by
\be\label{vsd}[\sd_{a,\dm I_d}v](k)=\dm^{d}\sum_{z\in\dZ}v(z)a(k-\dm z),\quad\forall k\in\dZ,\quad v\in\dlrs{0}{s}{r}.\ee
For simplicity, we only work with the case $s=1$, and the general case can be handled similarly. Given an initial data $v_0\in\dlrs{0}{1}{r}$, a \emph{vector $\dm I_d$-subdivision scheme} is implemented by applying the subdivision operator $\sd_{a,\dm I_d}$ on $v_0$ iteratively and generates a sequence $\{\sd_{a,\dm I_d}^nv_0\}_{n=1}^\infty$ of data. A subdivision scheme is closely related to a cascade algorithm. Define the \emph{Kronecker function $\td$} by
\be\label{td}\td(k)=\begin{cases}1, &k=0,\\
	0, &k\in\dR\setminus\{0\},\end{cases}\ee
and for simplicity, we also use $\td$ to denote its restriction $\td|_{\dZ}$. Define the \emph{semi-convolution} of a filter $v\in\dlrs{0}{s}{r}$ and an $r\times t$ matrix-valued function $f$ by
$$[v*f](x):=\sum_{k\in\dZ}v(k)f(x-k),\quad\forall x\in\dR.$$
It is easy to observe that
$$R_{a,\dm I_d}\psi(x)=[(\sd_{a,\dm I_d}(\td I_r))*\psi](\dm x)=\sum_{k\in\dZ}(\sd_{a,\dm I_d}(\td I_r))(k)\psi(\dm x-k),\quad\forall x\in\dR,$$
and more generally,
\be\label{ref:sd}R_{a,\dm I_d}^n\psi(x)=[(\sd_{a,\dm I_d}(\td I_r))*\psi](\dm^n x)=\sum_{k\in\dZ}(\sd_{a,\dm I_d}^n(\td I_r))(k)\psi(\dm^n x-k),\quad\forall n\in\N,\quad x\in\dR.\ee
From \er{ref:sd}, it is natural to choose the subdivision operator $\sd_{a,\dm I_d}$ to serve as a discrete version of the cascade operator $R_{a,\dm I_d}$.

\subsection{Convergence of Cascade Algorithms and Subdivision Schemes}

To use a vector cascade algorithm to compute a refinable vector function and its derivatives, we require the algorithm to converge in some space of smooth functions. For $m\in\N_0$ , denote $C^m(\dR)$ the space of \emph{$m$-times continuously differentiable functions}. Let $\dm \in\N\setminus\{1\}$, $r\in\N$ and $a\in \dlrs{0}{r}{r}$, we say that the vector $\dm I_d$-cascade algorithm with the filter $a$ is \emph{$C^m$-convergent} if for every suitably chosen (we will explain the requirements of an initial vector function in Section~\ref{sec:disc}) compactly supported initial vector function $\psi\in (C^m(\dR))^{r\times 1}$, the sequence $\{R_{a,\dm I_d}^n\psi\}_{n=1}^\infty$ converges to some compactly supported vector function $\psi_\infty\in (C^m(\dR))^{r\times 1}$, that is, 
$$\|R_{a,\dm I_d}^n\psi-\psi_\infty\|_{(C^m(\dR))^{r\times1}}:=\sum_{j=1}^r\sum_{\mu\in\cup_{q=0}^m\N_{0,q}^d}\left\|\partial^\mu\left[R_{a,\dm I_d}^n\psi-\psi_\infty\right]_j\right\|_{\dLp{\infty}}\to 0,$$
as $n\to\infty$, where
$$\N_{0,q}^d:=\{\nu=(\nu_1,\dots,\nu_d)\in\N_0^d:\, |\nu|=\nu_1+\dots+\nu_d=q\},\quad\forall q\in\N_0,$$
and $\partial^\mu:=\partial_1^{\mu_1}\dots\partial_d^{\mu_d}$ is the differential operator for all $\mu\in\N_0^d$. As shown in \cite[Theorem 4.3]{han03}, if the initial vector function $\psi$ satisfies some admissibility conditions (see Section~\ref{sec:disc}), then the limit function $\psi_\infty$ must be an $\dm I_d$-refinable vector function $\phi$ that satisfies \er{ref:eq}. In such cases, we have $\phi\in (C^m(\dR))^{r\times 1}$, and the cascade algorithm with the filter $a$ gives us a way to compute $\phi$ and its derivatives.

From \er{ref:sd}, we see that the vector cascade operator $R_{a,\dm I_d}$ and the vector subdivision operator $S_{a,\dm I_d}$ have intrinsic connections. Hence, it is reasonable to guess that the $C^m$-convergence of a vector cascade algorithm can be described using an analog of ''$C^m$-convergence" of its associated vector subdivision scheme. If $\{R_{a,\dm I_d}^n\psi\}_{n=1}^\infty$ converges to $\phi$ in $(C^m(\dR))^{r\times 1}$, then $\{\partial^\mu[R_{a,\dm I_d}^n\psi]\}_{n=1}^\infty$, or equivalently (using \er{ref:sd}), $\{[(\sd_{a,\dm I_d}^n(\td I_r))*(\dm^{|\mu|n}\partial^\mu\psi)](\dm^n\cdot)\}_{n=1}^\infty$ converges uniformly to $\partial^\mu\phi$ for all $\mu\in\cup_{q=0}^m\N_{0,q}^d$. Therefore, it is not surprising to guess that one suitable discrete version of describing the $C^m$-convergence of $\{R_{a,\dm I_d}^n\psi\}_{n=1}^\infty$ should be the following:
\be\label{defn:sd:0}\begin{aligned}&\text{For every }\mu\in\cup_{q=0}^m\N_{0,q}^d,\, \{[\sd_{a,\dm I_d}^n(\td I_r)]*u\}_{n=1}^\infty\text{  ``converges" to }\partial^\mu\phi\\
	&\text{ in some sense for some suitably chosen }u\in\dlrs{0}{r}{1},\end{aligned}\ee
where the \emph{convolution} of two filters $v\in\dlrs{0}{s}{r}$ and $w\in\dlrs{0}{r}{t}$ is defined via
$$[v*w](k):=\sum_{z\in\dZ}v(z)w(k-z)=\sum_{z\in\dZ}v(k-z)w(z),\quad\forall k\in\dZ.$$
For the scalar case $r=1$, the filter $u$ in \er{defn:sd:0} is obtained using a difference operator. For $h\in\dZ$, define the \emph{backward difference operator} $\nabla_h$ by:
\[
\nabla_h v(k)=v(k)-v(k-h),\quad\forall v\in \dlrs{0}{s}{r},\,k\in\dZ.
\]
For every $\mu=(\mu_1,\dots,\mu_d)\in\N_0^d$, we define
\[
\nabla^\mu:=\nabla^{\mu_1}_{\pe_1^d}\dots\nabla^{\mu_d}_{\pe_d^d},
\]
where $\pe_j^d$ is the $j$-th coordinate vector in $\dR$ for all $j=1,\dots,d$. For any $\mu\in\N_0^d$ and $v\in\dlrs{0}{s}{r}$, observe that $\nabla^\mu v =[\nabla^\mu \td]*v$ and it is straightforward to check that
\[
\wh{\nabla^\mu v}(\xi)=\wh{\nabla^\mu\td}(\xi)\wh{v}(\xi)=(1-e^{-i\xi_1})^{\mu_1}
\dots(1-e^{-i\xi_d})^{\mu_d}\wh{v}(\xi),
\quad\forall\xi=(\xi_1,\dots,\xi_d)\in\dR.
\]
We then have the following definition of a $C^m$-convergent scalar subdivision scheme (see e.g. \cite[Section 7.3.1]{hanbook}).
\begin{definition}\label{def:svd:0}[$C^m$-Convergence of a scalar subdivision scheme] Let $\dm\in\N\setminus\{0\}$, $m\in\N_0$ and $a\in\dlp{0}$ be such that $\wh{a}(0)=1$. If for every initial data $v\in\dlp{0}$, there exists $\eta_v\in C^m(\dR)$ such that
	\be\label{cm:conv:1}\lim_{n\to\infty}\sup_{k\in\dZ}\left|\dm^{|\mu|n}[\nabla^\mu\sd_{a,\dm}^nv](k)-\partial^\mu\eta_v(\dm^{-n}k)\right|=0,\quad\forall  \mu\in\cup_{q=0}^m\N_{0,q}^d, \ee
	then we say that the $\dm I_d$-subdivision scheme that employs the mask $a$ is $C^m$-convergent. In particular, if $v=\td$, then the limit function $\eta_{\td}=\phi$ must be the unique compactly supported function $\phi$ that satisfies $\wh{\phi}(\dm\xi)=\wh{a}(\xi)\wh{\phi}(\xi)$ and $\wh{\phi}(0)=1$. 
\end{definition}
Note that \er{cm:conv:1} with $v=\td$ can be equivalently written as
$$\lim_{n\to\infty}\sup_{k\in\dZ}\left|[(\sd_{a,\dm}^n\td)*(u_{\mu})](k)-\partial^\mu\phi(\dm^{-n}k)\right|=0,\quad\forall\mu\in\cup_{q=0}^m\N_{0,q}^d,$$
where $u_{\mu}:=\dm^{|\mu|n}\nabla^\mu\td$. Such a choice of $u_{\mu}$ is not that surprising since it is natural to use the divided difference operator $\nabla^\mu$ as a discrete version of the partial differential operator $\partial^\mu$ to measure the smoothness of functions. Furthermore, the scaling factor $\dm^{|\mu|n}$ is required to make $\dm^{|\mu|n}\nabla^\mu\sd_{a,\dm}^n\td\sim\partial^\mu\phi(\dm^{-n}\cdot)$ as $n$ gets large. With Definition~\ref{def:svd:0}, we have the following theorem on the equivalence between the convergences of a scalar subdivision scheme and a scalar cascade algorithm.

\begin{theorem}\label{thm:ssd:scs}[\cite[Theorem 7.3.1]{hanbook}] Let $\dm\in\N\setminus\{1\}$ and  $m\in\N_0$. Let $a\in\dlp{0}$ be such that $\wh{a}(0)=1$ and $\phi$ be the unique compactly supported distribution that satisfies $\wh{\phi}(\dm\xi)=\wh{a}(\xi)\wh{\phi}(\xi)$ and $\wh{\phi}(0)=1$. Then the following statements are equivalent:
	\begin{enumerate}
		
		\item[(1)] $\phi\in C^m(\dR)$ and $\{R_{a,\dm}^nf\}_{n=1}^\infty$ converges to $\phi$ in $C^m(\dR)$ for every $f\in C^m(\dR)$ that is compactly supported and satisfies
		$$\wh{f}(0)=1,\quad \partial^\mu\wh{f}(2\pi k)=0,\quad\forall \mu\in\cup_{q=0}^m\N_{0,q}^d,\quad k\in\dZ\setminus\{0\}.$$
		
		\item[(2)] The $\dm I_d$-subdivision scheme that employs the mask $a$ is $C^m$-convergent as in Definition~\ref{def:svd:0}.
	\end{enumerate}
	
\end{theorem}

Theorem~\ref{thm:ssd:scs} links the $C^m$-convergence of scalar cascade algorithms and subdivision schemes. Then, it is natural to ask whether we have similar results for the case $r>1$ on vector cascade algorithms and subdivision schemes. Unfortunately, so far, to the best of our knowledge, there is no notion of a so-called ``$C^m$-convergent vector subdivision scheme" that works for a general matrix-valued filter $a$ and a dilation matrix $\dm I_d$. For the case $d=1$ and $\dm =2$, a notion of a $C^m$-convergent vector subdivision scheme has been introduced and investigated in \cite{han23-2}. Motivated by the work of \cite{han23-2} and the importance of multivariate vector subdivision schemes in applications, our main goal of this paper is to develop a complete theory on the convergence of vector subdivision schemes that work for any matrix-valued filters in arbitrary dimensions.

\subsection{Our Contributions and Paper Structure}\label{ssec:con}

We have some clues from previous discussions to define the convergence of a vector subdivision scheme for the case $r>1$. For convenience, suppose $r,s\in\N$ and $(X,\|\cdot\|_X)$ is a normed $\C$-vector space, then for $f:=(f_{jk})_{1\le j\le r,\,1\le k\le s}\in X^{r\times s}$, define
$$\|f\|_{X^{r\times s}}:=\sum_{j=1}^r\sum_{k=1}^s\|f_{jk}\|_X.$$
Let $\dm\in\N\setminus\{1\}$, $r\in\N$, $m\in\N_0$ and $a\in\dlrs{0}{r}{r}$. For every initial data $v\in\dlrs{0}{1}{r}$, the vector subdivision scheme that uses the filter $a$ generates a sequence $\{\sd_{a,\dm I_d}^nv\}_{n=1}^\infty$. Motivated by the relations between the scalar cascade and subdivision operators and Definition~\ref{def:svd:0}, we guess that the definition of a``$C^m$-convergent vector subdivision scheme" takes the following form:
\begin{itemize}
	
	\item For every initial data $v\in\dlrs{0}{1}{r}$, there exists $\eta_v\in C^m(\dR)$ such that 
	\be\label{guess:vsd}\lim_{n\to\infty}\|(\sd_{a,\dm I_d}^nv)*(\dm^{\tau n}u_\mu)-\partial^\mu\eta_v(\dm^{-n}\cdot)\|_{\dlrs{\infty}{r}{1}}=0,\quad\forall \mu\in\cup_{q=0}^m\N_{0,q}^d,\ee
	where $u_\mu\in\dlrs{0}{r}{1}$ is suitably chosen and plays the role of $\nabla^\mu\td$ in a $C^m$-convergent scalar subdivision scheme, and $\dm^{\tau n}$ is a scaling factor for some $\tau\ge 0$. Moreover, the base function must be a refinable vector function $\phi=[\phi_1,\dots,\phi_d]^\tp$ that satisfies $\wh{\phi}(\dm\xi)=\wh{a}(\xi)\wh{\phi}(\xi)$. That is, if $v=\td [\pe_l^r]^\tp$ for some $l\in\{1,\dots,r\}$, then we must have $\eta=c\phi_l$ for some $c\in\C\setminus\{0\}$.\end{itemize}
We will see that the above guess is reasonable. Later, we will prove in  Theorem~\ref{thm:def:vsd} that there is not too much freedom to choose $u_\mu$ and $\tau$. Consequently, there is only one meaningful way to define a``$C^m$-convergent vector subdivision scheme. To do this, we need some notations. Let $m\in\N_0$, $\xi_0\in\dR$ and $f,g$ be two matrices (of the same size) of smooth functions. By 
$$f(\xi)=g(\xi)+\bo(\|\xi-\xi_0)\|^{m}),\quad \xi\to\xi_0,$$
we mean
$$\partial^\mu f(\xi_0)=\partial^\mu g(\xi_0),\quad\forall \mu\in\cup_{q=0}^{m-1}\N_{0,q}^d.$$
For every $\alpha\in\N_0^d$, $\vgu\in\dlrs{0}{1}{r}$ and $u\in\dlrs{0}{r}{1}$, define
\be\label{bvu} \beta_{\vgu,u,\mu}:=\frac{\partial^\mu[\wh{\vgu}\wh{u}](0)}{i^{|\mu|}\mu!},\ee
and define the linear space $\mom_{\vgu,\mu}$ as 
\be\label{mom}\mom_{\vgu,\mu}:=\{u\in\dlrs{0}{r}{1}:\, \wh{\vgu}(\xi)\wh{u}(\xi)=\beta_{\vgu,u,\mu}(i\xi)^\mu+\bo(\|\xi\|^{|\mu|+1}),\,\xi\to 0\}.\ee
Now, we are ready to present our first main contribution, which is the following definition of a $C^m$-convergent vector subdivision scheme with a general matrix-valued filter $a\in\dlrs{0}{r}{r}$.

\begin{definition}\label{def:vsd}Let $\dm\in\N\setminus\{1\}$, $r\in\N$ and $m\in\N_0$. Let $a\in\dlrs{0}{r}{r}$ be a finitely supported filter that has order $m+1$ sum rules with respect to $\dm I_d$, that is, there exists an order $m+1$ matching filter $\vgu_a\in\dlrs{0}{1}{r}$ of $a$ that satisfies
	\be\label{sr}\wh{\vgu_a}(0)\ne0,\quad \wh{\vgu_a}(\dm\xi)\wh{a}(\xi+2\pi\omega)=\td(\omega)\wh{\vgu_a}(\xi)+\bo(\|\xi\|^{m+1}),\quad\xi\to 0,\quad \forall \omega\in\Omega_{\dm I_d},\ee
	where $\td$ is defined via \er{td} and
	\be\label{om:dm}\Omega_{\dm I_d}:=[\dm^{-1} \dZ]\cap[0,1)^d.\ee
	For any given $u\in\dlrs{0}{r}{1}$ and initial data $v\in\dlrs{0}{1}{r}$, a vector $\dm I_d$-subdivision scheme with the filter $a$ generates a sequence $\{[\sd_{a,\dm I_d}^nv]*u\}_{n=1}^\infty$. The vector $\dm I_d$-subdivision scheme with the filter $a$ is \emph{$C^m$-convergent} if for every initial data $v\in\dlrs{0}{1}{r}$, there exists $\eta_v\in C^m(\dR)$ such that 
	\be\label{conv:vsd}\lim_{n\to\infty}\left\|\dm^{|\mu|n}(\sd_{a,\dm I_d}^nv)*u-\beta_{\vgu_a,u,\mu}\partial^\mu\eta_v(\dm^{-n}\cdot)]\right\|_{\dlp{\infty}}=0,\quad \forall u\in\mom_{\vgu_a,\mu},\, \mu\in\cup_{q=0}^m\N_{0,q}^d,\ee
	where $\beta_{\vgu_a,u,\mu}$ is defined as \er{bvu} and $\mom_{\vgu_a,\mu}$ is defined as \er{mom} with $\vgu=\vgu_a$.
\end{definition}

Next, our second main contribution is the following Theorem, which provides a comprehensive characterization of the convergence and smoothness of a vector subdivision scheme.

\begin{theorem}\label{thm:conv:vsd}Let $\dm\setminus\{1\}$, $r\in\N$ and $m\in\N_0$. Let $a\in\dlrs{0}{r}{r}$ be a filter with an order $m+1$ matching filter $\vgu_a\in\dlrs{0}{1}{r}$ that satisfies \er{sr}. If $a$ satisfies
	\be\label{eigen:a}1\text{ is a simple eigenvalue of }\wh{a}(0)\text{ and all its other eigenvalues are smaller than }\dm^{-m}\text{ in modulus},\ee 
	it is well-known that there exists a unique $r\times 1$ vector $\phi$ of compactly supported distributions that satisfies \be\label{stan:vref}\wh{\phi}(\dm\xi)=\wh{a}(\xi)\wh{\phi}(\xi)\text{ and }\wh{\vgu_a}(0)\wh{\phi}(0)=1.\ee 
	Then the following statements are equivalent:
	
	\begin{enumerate}
		
		\item[(1)]$\sm_\infty(a,\dm I_d)>m$, where $\sm_\infty(a,\dm I_d)$ is the \emph{$L_\infty$-smoothness exponent of the filter $a$ with respect to $\dm I_d$}, and is defined in Section~\ref{sec:disc}. Consequently, the $\dm I_d$ vector cascade algorithm with the filter $a$ is $C^m$-convergent (i.e., item (7) of Theorem~\ref{thm:vcs} holds).

		\item[(2)] $\phi\in (C^m(\dR))^{r\times 1}$ and
		\be\label{conv:vsd:phi}\lim_{n\to\infty}\left\|\dm^{|\mu|n}(\sd_{a,\dm I_d}^n(\td I_r))*u-\beta_{\vgu_a,u,\mu}\partial^\mu\phi(\dm^{-n}\cdot)]\right\|_{\dlrs{\infty}{r}{1}}=0,\ee
		for all $u\in \mom_{\vgu_a,\mu},$ $\mu\in\cup_{q=0}^m\N_{0,q}^d$.
		
		\item[(3)] $\phi\in (C^m(\dR))^{r\times 1}$ and \er{conv:vsd:phi} holds for all $u\in G_{\vgu_a,\mu},$ $\mu\in\cup_{q=0}^m\N_{0,q}^d$, where $G_{\vgu_a,\mu}$ is any set that generates $\mom_{\vgu_a,\mu}$, that is, every $w\in\mom_{\vgu_a,\mu}$ can be written as $w=w_1*c_1+\dots+w_t*c_t$ for some $w_1,\dots,w_t\in G_{\vgu_a,\mu}$ and $c_1,\dots,c_t\in\dlp{0}$.

		\item[(4)] $\phi\in (C^m(\dR))^{r\times 1}$ and there exists $\mu\in\N_{0,m}^d$ such that \er{conv:vsd:phi} holds for all $u\in G_{\vgu_a,\mu}$, where $G_{\vgu_a,\mu}$ is any set that generates $\mom_{\vgu_a,\mu}$.

		\item[(5)] The vector $\dm I_d$-subdivision scheme that uses the filter $a$ is $C^m$-convergent as in Definition~\ref{def:vsd}.

	\end{enumerate}
	
	If $\phi\in(C^m(\dR))^{r\times 1}$ and $\phi$ has stable integer shifts, that is, $\spa\{\wh{\phi}(\xi+2\pi k):\,k\in\dZ\}=\C^{r\times 1}$ for all $\xi\in\dR$, then all items (1)-(5) must hold.
	
\end{theorem}

The paper is organized as follows:

\begin{itemize}
	\item Section~\ref{sec:disc} provides some necessary backgrounds on vector cascade algorithms. First, we will explain the role of a matching filter $\vgu_a$ of the filter $a$ and how it relates to the initial function of a vector cascade algorithm. Next, we will give a detailed discussion on the $L_p$-smoothness exponent $\sm_p(a,\dm I_d)$ of the filter $p$, which is the key to characterizing the convergence and smoothness of a vector cascade algorithm. Notably, we will prove Theorem~\ref{thm:sm} on the relation between $\sm_p(a,\dm I_d)$ and the matching filters of $a$, which, later, plays a key role in analyzing the convergence rate of a vector subdivision scheme.

	\item Section~\ref{sec:def} provides a detailed justification of Definition~\ref{def:vsd}. Specifically, we will prove Theorem~\ref{thm:def:vsd}, which tells us what conditions the filter $u_\mu$ and the scaling exponent $\tau$ in \er{guess:vsd} need to satisfy. As a consequence, Theorem~\ref{thm:def:vsd} explains why Definition~\ref{def:vsd} is the only meaningful way to define a vector subdivision scheme and its convergence.

	\item Section~\ref{sec:main} performs a detailed analysis of our newly defined vector subdivision scheme and its convergence. We first prove the main result Theorem~\ref{thm:conv:vsd}. Then, we will prove in Theorem~\ref{thm:cr} how fast a vector subdivision scheme converges. Next, we will discuss transformations of vector subdivision schemes. Specifically, we show that any vector subdivision scheme can be transformed into a scalar-type subdivision scheme. Last but not least, we will apply our theory of vector subdivision schemes to discuss the convergence of Lagrange and Hermite subdivision schemes.

	\item Section~\ref{sec:exmp} provides several examples of  two-dimensional vector $2I_2$-subdivision schemes to illustrate our theory.

\end{itemize}

\section{Convergence of Vector Cascade Algorithms}\label{sec:disc}

Before we justify Definition~\ref{def:vsd} and prove the main result Theorem~\ref{thm:conv:vsd}, let us discuss the convergence of vector cascade algorithms in more detail. Since vector cascade algorithms and subdivision schemes are closely related, several properties of vector cascade algorithms are powerful tools for the analysis of vector subdivision schemes. Here, we first briefly discuss the role of a matching filter $\vgu_a$ and how it relates to the initial function of a vector cascade algorithm. Then, we will perform detailed discussions on the critical quantity $\sm_p(a,\dm I_d)$, which plays a key role in the convergence of vector cascade algorithms and subdivision schemes.

\subsection{The Matching Filter $\vgu_a$ and the Conditions for an Initial Function}

Let $\dm\in\N\setminus\{1\}$,  $r\in\N$ and $m\in\N_0$. Suppose $a\in\dlrs{0}{r}{r}$ is a filter that has order $m+1$ sum rules with a matching filter $\vgu_a\in\dlrs{0}{1}{r}$ that satisfies \er{sr}. Using the Leibniz differentiation rule, \er{sr} yields
\be\label{sr:diff}\partial^\mu \wh{\vgu_a}(0)=\sum_{\nu\le\mu}\binom{\mu}{\nu}\dm^{|\nu|}\partial^{\nu}\wh{\vgu_a}(0)\partial^{\mu-\nu}\wh{a}(0),\quad\forall \mu\in\cup_{q=0}^m\N_{0,q}^d,\ee
where for $\nu=(\nu_1,\dots,\nu_d),\mu=(\mu_1,\dots,\mu_d)\in\N_0^d$, $\nu\le\mu$ means $\nu_t\le \mu_t$ for all $t=1,\dots,d.$ Note that an order $m+1$ matching filter is fully determined by the values of the derivatives $\partial^{\mu}\wh{\vgu}(0)$ for all $\mu\in\cup_{q=0}^m\N_{0,q}^d$. Given a filter $a$, we simply solve the system of linear equations in \er{sr:diff} to get the values of $\partial^{\mu}\wh{\vgu}(0)$ for all $\mu\in\cup_{q=0}^m\N_{0,q}^d$. If in addition, the filter $a$ satisfies the mild condition in \er{eigen:a}, then $I_r-\dm^{j}\wh{a}(0)$ is invertible whenever $j\in\{1,\dots,m\}$. Hence, by \er{sr:diff}, up to a multiplicative constant, the set $\{\partial^{\mu}\wh{\vgu}(0):\,\mu\in\cup_{q=0}^m\N_{0,q}^d\}$ is uniquely determined through the following recursive formulas:
\be\label{sr:diff:0}\wh{\vgu_a}(0)\wh{a}(0)=\wh{\vgu_a}(0),\quad\partial^{\mu}\wh{\vgu_a}(0)=\sum_{\nu<\mu}\binom{\mu}{\nu}\partial^{\nu}\wh{\vgu_a}(0)\partial^{\mu-\nu}\wh{a}(0)[I_r-\dm^{|\mu|}\wh{a}(0)]^{-1},\quad \mu\in\cup_{q=1}^m\N_{0,q}^d.\ee
The condition \er{eigen:a} also yields that there exists a unique $r\times 1$ vector $\phi$ of compactly supported distributions such that \er{stan:vref} holds. We call the unique vector $\phi$ of compactly supported distributions that satisfies \er{stan:vref} the \emph{standard $\dm I_d$-refinable function vector of the filter $a$ associated with $\vgu_a$.}\\

\vsp{0.2}

Next, we discuss the relation between $\vgu_a$ and any initial function $\psi$ of a vector cascade algorithm. To study the smoothness of refinable vector functions, it is essential to investigate the convergence of vector cascade algorithms with a smooth initial function. For $m\in\N_0$ and $p\in[1,\infty]$, denote $W_p^m(\dR)$ the \emph{$L_p$-Sobolev space} which consists of all functions $f$ such that $\partial^\mu f\in\dLp{p}$ for all $\mu\in\cup_{q=0}^m\N_{0,q}^d$, endowed with the following norm
$$\|f\|_{W_p^m(\dR)}:=\sum_{\mu\in\cup_{q=0}^m\N_{0,q}^d}\|\partial^\mu f\|_{\dLp{p}}.$$ 
The following result gives the essential conditions on an initial function $\psi$ that generates a convergent sequence $\{R_{a,\dm I_d}^n\psi\}_{n=1}^\infty$ in the Sobolev space $(W_p^m(\dR))^{r\times 1}$.

\begin{prop}\label{prop:in}[\cite[Propositions 3.1 and 3.2]{han03}]Let $\dm \in\N\setminus\{1\}$, $r\in\N$, $m\in\N_0$ and $p\in[1,\infty]$. Let $a\in\dlrs{0}{r}{r}$ be a filter that has order $m+1$ sum rules with a matching filter $\vgu_a\in\dlrs{0}{1}{r}$ that satisfies \er{sr}. Then the following statements hold:
	
	\begin{enumerate}
		\item[(1)] Let $\psi\in(W_p^m(\dR))^{r\times 1}$ be a vector of compactly supported functions such that $\spa\{\wh{\psi}(\xi+2\pi k):\,k\in\dZ\}=\C^{r\times 1}$. If $\lim_{n\to\infty}\|R_{a,\dm I_d}^n\psi-\psi_\infty\|_{(W_p^m(\dR))^{r\times 1}}=0$ for some $\psi_\infty\in (W_p^m(\dR))^{r\times 1}\setminus\{0\}$, then \er{eigen:a} must hold. Consequently, up to a multiplicative constant, the set $\{\partial^\mu\wh{\vgu_a}(0):\,\mu\in\cup_{q=0}^m\N_{0,q}^d\}$ is uniquely determined through \er{sr:diff:0}.
		
		\item[(2)]For any compactly supported vector $\psi\in(W_p^m(\dR))^{r\times 1}$, if $\{R_{a,\dm I_d}^n\psi\}_{n=1}^\infty$ is convergent in $(W_p^m(\dR))^{r\times 1}$ and $\lim_{n\to\infty}\wh{\vgu_a}(0)\wh{R_{a,\dm I_d}^n\psi}(0)=1$, then $\psi$ must satisfy the \emph{order $m+1$ admissibility condition}, that is,
		\be\label{psi:adm}\wh{\vgu_a}(0)\wh{\psi}(0)=1,\quad  \wh{\vgu_a}(\xi)\wh{\psi}(\xi+2\pi k)=\bo(\|\xi\|^{m+1}),\quad\xi\to 0,\quad\forall k\in\dZ\setminus\{0\}.\ee
		If, in addition, \er{eigen:a} holds, then there exists a unique $r\times 1$ vector $\phi$ of compactly supported distributions that satisfies \er{stan:vref}. When $\phi\in(W_p^m(\dR))^{r\times 1}$, it must satisfy the following admissibility condition:
		\be\label{phi:adm}\wh{\vgu_a}(\xi)\wh{\phi}(\xi+2\pi k)=\td(k)+\bo(\|\xi\|^{m+1}),\quad\xi\to 0,\quad\forall k\in\dZ.\ee
	\end{enumerate}
\end{prop}

\subsection{Convergence of Vector Cascade Algorithms and the Quantity $\sm_p(a,\dm I_d)$}
To characterize the convergence of vector cascade algorithms, we need the technical quantity $\sm_p(a,\dm I_d)$, known as the \emph{$L_p$-smoothness exponent of the filter $a$ with respect to $\dm I_d$}. Let $m\in\N_0$, $p\in[1,\infty]$ and $y\in \dlrs{0}{1}{r}$ be such that $\wh{y}(0)\ne 0$, define
\be\label{vmy}V_{-1,y}:=\dlrs{0}{r}{1},\quad V_{m,y}:=\{w\in\dlrs{0}{r}{1}:\, \wh{y}(\xi)\wh{w}(\xi)=\bo(\|\xi\|^{m+1}),\,\xi\to 0\},\ee
and define 
\be\label{rho:a:0}\rho_{m}(a,\dm I_d,y)_p:=\sup\left\{\limsup_{n\to\infty}\|[\sd_{a,\dm I_d}^n(\td I_r)]*w\|_{\dlrs{p}{r}{1}}^{\frac{1}{n}}:\, w\in V_{m-1,y}\right\}.\ee
As $\wh{y}(0)=[\wh{y_1}(0),\dots,\wh{y_r}(0)]\ne 0$, without loss of generality, assume $\wh{y_1}(0)\ne 0$. Then for every $j=2,\dots,r$, choose $b_j\in\dlp{0}$ such that 
$$\wh{b_j}(\xi)=\frac{\wh{y_j}(\xi)}{\wh{y_1}(\xi)}+\bo(\|\xi\|^{m}),\quad\xi\to 0,\quad j=2,\dots,r,$$
and define
$$B_{m-1,y}:=\{\nabla^{\mu}\td \pe_1^r:\, \mu\in\N_{0,m}^d\}\cup\{b_j\pe^r_1+\td \pe^r_j:\, j=2,\dots,r\},$$
we can see that $B_{m-1,y}$ generates $V_{m-1,y}$, that is, every $w\in V_{m-1,y}$ can be written as $w=w_1*c_1+\dots+w_s*c_s$ for some $w_1,\dots,w_s\in B_{m-1,y}$ and $c_1,\dots,c_s\in\dlp{0}$. Therefore, one can further deduce that
\be\label{rho:a:1}\rho_{m}(a,\dm I_d,y)_p:=\max\left\{\limsup_{n\to\infty}\|[\sd_{a,\dm I_d}^n(\td I_r)]*w\|_{\dlrs{p}{r}{1}}^{\frac{1}{n}}:\, w\in B_{m-1,y}\right\}.\ee
Now define
\be\label{sr:a}\begin{aligned}m_a:=\sr(a,\dm I_d):=\sup\{&m\in\N_0:\,\text{there exists }y\in\dlrs{0}{1}{r}\text{ such that }\wh{y}(0)\ne 0\\
	&\text{ and } \wh{y}(\dm\xi)\wh{a}(\xi+2\pi\omega)=\td(\omega)\wh{y}(\xi)+\bo(\|\xi\|^m)\text{ for all }\omega\in\Omega_{\dm I_d}\},\end{aligned}\ee
and 
\be\label{rho:a:00}\begin{aligned}\rho(a,\dm I_d)_p:=\inf\{&\rho_m(a,\dm I_d,y):\,m\in\{0,\dots,m_a\}\text{ and }y\in\dlrs{0}{1}{r}\text{ satisfy}\\
	&\wh{y}(0)\ne 0\text{ and } \wh{y}(\dm\xi)\wh{a}(\xi+2\pi\omega)=\td(\omega)\wh{y}(\xi)+\bo(\|\xi\|^m)\text{ for all }\omega\in\Omega_{\dm I_d}\}.\end{aligned}\ee
We define the \emph{$L_p$-smoothness exponent of $a$ with respect to $\dm I_d$} by
\be\label{sm:a}\sm_p(a,\dm I_d):=\frac{d}{p}-\log_{\dm}\rho(a,\dm I_d)_p.\ee
The quantity $\sm_p(a,\dm I_d)$ plays a key role in analyzing the convergence of the vector cascade algorithm with the filter $a$ in the $L_p$-Sobolev space. Here, we recall the following important result from \cite{han03}, which gives a detailed characterization of the convergence of vector cascade algorithms in Sobolev spaces.

\begin{theorem}\label{thm:vcs}[\cite[Theorem 4.3]{han03}]. Let $\dm\in\N\setminus\{1\}$, $r\in\N$, $m\in\N_0$ and $p\in[1,\infty]$. Let $a\in\dlrs{0}{r}{r}$ and $\vgu_a\in\dlrs{0}{1}{r}$ be such that \er{sr} holds and $\wh{\vgu_a}(0)\ne 0$. Then the following statements are equivalent:
	
	\begin{enumerate}
		\item[(1)] The vector $\dm I_d$-cascade algorithm with the filter $a$ converges in $W_p^m(\dR)$, that is, $\{R_{a,\dm I_d}^n\psi\}_{n=1}^\infty$ is a Cauchy sequence in $(W_p^m(\dR))^{r\times 1}$ for every $\psi\in (W_p^m(\dR))^{r\times 1}$ that is compactly supported and satisfies order $m+1$ admissibility condition as in \er{psi:adm}.

		\item[(2)]For some function $\psi\in (W_p^m(\dR))^{r\times 1}$ (require $\psi\in (C^m(\dR))^{r\times 1}$ if $p=\infty$) that satisfies the admissibility condition in \er{psi:adm} and $\spa\{\wh{\psi}(\xi+2\pi k):\, k\in\dZ\}=\C^{r\times 1}$ for all $\xi\in\dR$, $\{R_{a,\dm I_d}^n\psi\}_{n=1}^\infty$ is a Cauchy sequence in $(W_p^m(\dR))^{r\times 1}$.
		
		\item[(3)]$\lim_{n\to\infty}\dm^{n(m-\frac{d}{p})}\|[\sd_{a,\dm I_d}^n(\td I_r)]*w\|_{\dlp{p}}=0$ for all $w\in V_{m,\vgu_a}$, where $V_{m,\vgu_a}$ is defined via \er{vmy} with $y=\vgu_a$.
		
		\item[(4)] $\rho_{m+1}(a,\dm I_d,\vgu_a)_p<\dm^{\frac{d}{p}-m}$.
		
		\item[(5)]$\sm_p(a,\dm I_d)>m$.
	\end{enumerate}
	If $\phi\in(W_p^m(\dR))^{r\times 1}$ (require $\phi\in (C^m(\dR))^{r\times1}$ if $p=\infty$) and $\phi$ has stable integer shifts, then all items (1)-(5) must hold. Moreover, any of the above items (1)-(5) implies 
	
	\begin{enumerate}
		\item[(6)]$\sr(a,\dm I_d)\ge\sm_p(a,\dm I_d)>m$.
		
		\item[(7)]The filter $a$ satisfies \er{eigen:a}. Consequently, there exists a unique $r\times 1$ vector $\phi$ of compactly supported functions that satisfies \er{stan:vref}. Furthermore, $\phi\in(W_p^m(\dR))^{r\times 1}$ (replace by $\phi\in(C^m(\dR))^{r\times 1}$ if $p=\infty$) and $\lim_{n\to\infty}\|R_{a,\dm I_d}^n \psi-\phi\|_{(W_p^m(\dR))^{r\times 1}}=0$ for every $\psi\in(W_p^m(\dR))^{r\times 1}$ that is compactly supported and satisfies the admissibility condition in \er{psi:adm}.
	\end{enumerate}
\end{theorem}

To compute $\sm_p(a,\dm I_d)$, it suffices to find $\rho(a,\dm I_d)_p$ as in \er{rho:a:00}, which seems to be difficult since $\rho(a,\dm I_d)_p$ is obtained by taking the infimum of $\rho_{m}(a,\dm I_d,y)_p$ over all $m\in\{0,\dots,m_a\}$ and all possible matching filters $y\in\dlrs{0}{1}{r}$ of $a$. Nevertheless, the following theorem tells us it suffices to pick a specific sum rule order $m$ and an associated matching filter to obtain $\rho(a,\dm I_d)_p$.

\begin{theorem}\label{thm:sm}Let $\dm\in\N\setminus\{1\}$, $r\in\N$ and $p\in[1,\infty]$. Let $a\in\dlrs{0}{r}{r}$ be such that $\sr(a,\dm I_d)=m_a$. Suppose $k_0\in\N_0$ is a non-negative integer with $\sm_p(a,\dm I_d)\le k_0\le m_a$ and $\vgu\in\dlrs{0}{1}{r}$ is an order $k_0$ matching filter of $a$ that satisfies
	\be\label{sr:vgu}\wh{\vgu}(0)\ne 0,\quad \wh{\vgu}(\dm\xi)\wh{a}(\xi+2\pi\omega)=\td(\omega)\wh{\vgu}(\xi)+\bo(\|\xi\|^{k_0}),\quad\xi\to 0,\quad \forall\omega\in\Omega_{\dm I_d}.\ee  Then
	\be\label{rho:vgu}\rho(a,\dm I_d)_p=\rho_{k_0}(a,\dm I_d,\vgu)_p.\ee
	Consequently, we have
	\be\label{sm:sr}\sm_p(a,\dm I_d)=\frac{d}{p}-\log_{\dm}\rho_{k_0}(a,\dm I_d,\vgu)_p.\ee
\end{theorem}

\bp Consider the following two cases:

\begin{itemize}
	
	\item Case 1: $\sm_p(a,\dm I_d)>0$. In this case, there exists $m\in\N_0$ such that $m<\sm_p(a,\dm I_d)\le m+1$. By Theorem~\ref{thm:vcs}, the filter $a$ must satisfy \er{eigen:a}. Let $\vgu_a\in\dlrs{0}{1}{s}$ be an order $m_a$ matching filter of $a$ that satisfies 
	\be\label{sr:vgu:a}\wh{\vgu_a}(0)\ne 0,\quad \wh{\vgu_a}(\dm\xi)\wh{a}(\xi)=\td(\omega)\wh{\vgu_a}(\xi)+\bo(\|\xi\|^{m_a}),\quad\xi\to 0,\quad\forall\omega\in\Omega_{\dm I_d}.\ee
	If $k_0\in\{0,\dots,m_a\}$ and $\vgu\in\dlrs{0}{1}{r}$ is an order $k_0$ matching filter of $a$ that satisfies \er{sr:vgu}, then by \er{eigen:a}, there must exist $c\in\C\setminus\{0\}$ such that $\wh{\vgu}(\xi)=c\wh{\vgu_a}(\xi)+\bo(\|\xi\|^{k_0})$ as $\xi\to 0$, which implies $V_{k_0-1,\vgu}=V_{k_0-1,\vgu_a}$ and $$\rho_{k_0}(a,\dm I_d,\vgu)=\rho_{k_0}(a,\dm I_d,\vgu_a).$$ 
	Hence, by the definition of $\rho(a,\dm I_d)_p$ in \er{rho:a:00}, we have
	\be\label{rho:a:01}\rho(a,\dm I_d)_p=\min_{k_0\in\{0,\dots,m_a\}}\{\rho_{k_0}(a,\dm I_d,\vgu_a)_p\}.\ee
	As pointed out after \cite[Lemma 4.1]{han03} (also see \cite[Theorem 3.1]{han03-2}), we have
	\be\label{rho:ko:ma}\rho_{k_0}(a,\dm I_d,\vgu_a)_p=\max\{\rho_{m_a}(a,\dm I_d,\vgu_a)_p,\,\dm^{\frac{d}{p}-k_0}\},\quad\forall k_0=0,\dots,m_a-1,\ee
	from which we see that
	$$\rho(a,\dm I_d)_p=\rho_{m_a}(a,\dm I_d,\vgu_a)_p.$$
	If $k_0\ge\sm_p(a,\dm I_d)$, then
	$$\rho_{m_a}(a,\dm I_d,\vgu_a)_p=\rho(a,\dm I_d)_p=\dm^{\frac{d}{p}-\sm_p(a,\dm I_d)}\ge\dm^{\frac{d}{p}-k_0}.$$
Consequently, we deduce that
	$$\rho_{k_0}(a,\dm I_d,\vgu)_p=\rho_{k_0}(a,\dm I_d,\vgu_a)_p=\rho_{m_a}(a,\dm I_d,\vgu_a)_p,\quad\forall k\in\N_0\text{ with }\sm_p(a,\dm I_d)\le k\le m_a.$$
	Therefore, \er{rho:vgu} must hold, and thus, \er{sm:sr} holds.

	\item Case 2: $\sm_p(a,\dm I_d)\le 0$. Let $\vgu_a\in\dlrs{0}{1}{r}$ be a filter that satisfies \er{sr:vgu:a}. By the definition of $\rho(a,\dm I_d)_p$ and $\dm >1$, we have
	$$\log_{\dm}\rho_{m_a}(a,\dm I_d,\vgu_a)_p\ge \log_{\dm}\rho(a,\dm I_d)_p=\frac{d}{p}-\sm_p(a,\dm I_d),$$
	which implies
	$$\rho_{m_a}(a,\dm I_d,\vgu_a)_p\geq \dm^{\frac{d}{p}-\sm_p(a,\dm I_d)}\ge \dm^{\frac{d}{p}-j},\quad\forall j\ge 0.$$
	Hence, it follows from \er{rho:ko:ma} that
	\be\label{rho:ko:ma:1}\rho_{k_0}(a,\dm I_d,\vgu_a)_p=\rho_{m_a}(a,\dm I_d,\vgu_a)_p,\quad\forall k_0=0,\dots,m_a.\ee
	For every $k_0\in\{0,\dots, m_a\}$, let $\vgu\in\dlrs{0}{1}{r}$ be an order $k_0$ matching filter of $a$ that satisfies \er{sr:vgu}. From \er{vmy}, we see that $V_{-1,\vgu}=V_{-1,\vgu_a}=\dlrs{0}{1}{r}$, which together with \er{rho:ko:ma} yield
	$$\rho_0(a,\dm I_d,\vgu_a)_p=\rho_0(a,\dm I_d,\vgu)_p=\max\{\rho_{k_0}(a,\dm I_d, \vgu)_p,\,\dm^{\frac{d}{p}}\}.$$
	By the assumption $\sm_p(a,\dm I_d)\le 0$, we have $\frac{d}{p}-\log_{\dm}\rho_{k_0}(a,\dm I_d,\vgu)_p\le 0$ and thus $\rho_{k_0}(a,\dm I_d,\vgu)_p\ge \dm^{\frac{d}{p}}$. Hence, 
	\be\label{sr:k0}\rho_0(a,\dm I_d,\vgu_a)_p=\rho_0(a,\dm I_d,\vgu)_p=\rho_{k_0}(a,\dm I_d, \vgu)_p,\quad \forall k_0=0,\dots,m_a,\quad \vgu\in\dlrs{0}{1}{r}\text{ satisfying \er{sr:vgu}}.\ee
	Therefore, we conclude from \er{rho:ko:ma:1}, \er{sr:k0} and the definition of $\rho(a,\dm I_d)_p$ that \er{rho:vgu} must hold and consequently, \er{sm:sr} also holds.
\end{itemize}\ep

Theorem~\ref{thm:sm} greatly simplifies the definition of $\sm_p(a,\dm I_d)$, but it is still not generally easy to compute this quantity. Indeed, for a general matrix filter $a\in\dlrs{0}{r}{r}$, there only exist methods to compute the $L_2$-smoothness exponent $\sm_2(a,\dm I_d)$, see, for instance, \cite[Theorem 7.1]{han03}, \cite[Theorem 2.4]{han03-2} and \cite[Theorem 1.1]{jj03}. For $p>2$, we can use the following inequalities to estimate $\sm_p(a,\dm I_d)$ from $\sm_2(a,\dm I_d)$ (see \cite{han03,han03-2}):
\be\label{sm:p2}\sm_2(a, \dm I_d)-d\left(\frac{1}{2}-\frac{1}{p}\right)\le\sm_p(a,\dm I_d)\le\sm_2(a,\dm I_d).\ee

\section{Justification of Definition~\ref{def:vsd}}\label{sec:def}

In this section, we provide a detailed justification of Definition~\ref{def:vsd}, that is, why using this definition to describe the convergence and smoothness of a vector subdivision scheme. As discussed at the beginning of Subsection~\ref{ssec:con}, a reasonable way to define a $C^m$-convergent subdivision scheme is through \er{guess:vsd} with a suitably chosen $u_\mu$ that plays the role of $\nabla^\mu\td$ in a $C^m$-convergent scalar subdivision scheme and a scaling factor $\dm^{\tau n}$. The following theorem tells us the requirements that $u_\mu$ and $\tau$ must satisfy to make the definition of a $C^m$-convergent vector subdivision scheme meaningful.

\begin{theorem}\label{thm:def:vsd}Let $\dm\in\N\setminus\{1\}$, $r\in\N$, $m\in\N_0$ and $u\in\dlrs{0}{r}{1}$. Let $a\in\dlrs{0}{r}{r}$ be a finitely supported filter and $\vgu_a\in\dlrs{0}{1}{r}$ be an order $m+1$ matching filter of $a$ that satisfies \er{sr}. Suppose there exist $\tau\in[0,\infty)$, $\mu\in\N_{0,m_0}^d$ for some $m_0\in\{0,\dots,m\}$ and $\eta\in (C^{m_0}(\dR))^{r\times 1}$ such that
	\be\label{conv:vsd:0}\lim_{n\to\infty}\|\dm^{\tau n}[\sd_{a,\dm I_d}^n(\td I_r)]*u-\partial^\mu \eta(\dm^{-n}\cdot)\|_{\dlrs{\infty}{r}{1}}=0.\ee
	Then the following statements hold:
	
	\begin{enumerate}
		
		\item[(1)] For all $\nu\in\cup_{q=0}^m\N_{0,q}^d$ that satisfy $|\nu|\le\tau$, we have
		\be\label{vgu:u:0}\partial^\nu[\wh{\vgu_a}\wh{u}](0)=\begin{cases}0, &\text{ if }|\nu|<\tau,\\
			\partial^\nu[\wh{\vgu_a}\wh{\partial^\mu\eta}](0), &\text{ if }|\nu|=\tau\in\N_0,\end{cases}.\ee
		
		\item[(2)]The function $\partial^\mu\eta$ in \er{conv:vsd:0} is compactly supported and is given by
		\be\label{ref:eta}\wh{\partial^\mu\eta}(\xi):=\lim_{n\to\infty}\left[\dm^{\tau n}\left(\prod_{j=1}^n\wh{a}(\dm^{-j}\xi)\right)\wh{u}(\dm^{-n}\xi)\right],\quad\xi\in\dR.\ee
		Moreover, $\partial^\mu \eta$ satisfies the refinable relation
		\be\label{ref:eta:0}\wh{\partial^\mu\eta}(\dm\xi)=\dm^{\tau }\wh{a}(\xi)\wh{\partial^\mu\eta}(\xi),\quad\xi\in\dR.\ee
		
		\item[(3)]  $\lim_{n\to\infty}\|\dm^{\tau n}[\sd_{a,\dm I_d}^n(\td I_r)*(u*w)-\wh{w}(0)\partial^\mu\eta(\dm^{-n}\cdot)]_{\dlrs{\infty}{r}{1}}=0$ for all $w\in\dlp{0}$.
		
	\end{enumerate}
	
	Suppose, in addition, $\tau\in[0,m]$, $\partial^\mu\eta$ is not identically zero and the filter $a$ satisfies \er{eigen:a}. Let $\phi$ be the unique $r\times 1$ vector of compactly supported distributions that satisfies $\wh{\phi}(\dm\xi)=\wh{a}(\xi)\wh{\phi}(\xi)$ and $\wh{\vgu_a}(0)\wh{\phi}(0)=1.$ Then 
	
	\begin{enumerate}		
		
		\item[(4)] $\tau\ge |\mu|=m_0$ is given by
		\be\label{tau}\begin{aligned}\tau=\inf\{N\in\N_0:\,& \partial^\nu[\wh{\partial^\mu\eta}](0)\ne 0\text{ for some }\nu\in\N_{0,N}^d\text{ and }\\
			&\partial^\alpha[\wh{\partial^\mu\eta}](0)=0\text{ for all }\alpha\in\cup_{q=0}^{N-1}\N_{0,q}^d\}.\end{aligned}\ee
		Particularly, $\tau=|\mu|$ if $\wh{\eta}(0)\ne 0$ and in this case, $\phi\in(C^{m_0}(\dR))^{r\times 1}$ and we must have $\eta=c\phi$ for some $c\in\C\setminus\{0\}$.
		
		\item[(5)]	The filter $u$ must satisfy
		\be\label{mom:vgu:u}\wh{\vgu_a}(\xi)\wh{u}(\xi)=\beta_{\vgu_a,u,\mu}(i\xi)^\mu+\bo(\|\xi\|^{m_0+1}),\quad\xi\to 0,\ee
		where $\beta_{\vgu_a,u,\alpha}:=\frac{\partial^\alpha[\wh{\vgu_a}\wh{u}](0)}{i^{m_0}\alpha!}$. Moreover, for every input data $v\in\dlrs{0}{1}{r}$, we have
		\be\label{conv:vsd:01}\lim_{n\to\infty}\left\|\dm^{\tau n}[\sd_{a,\dm I_d}^nv]*u-\partial^\mu[v*\eta](\dm^{-n}\cdot)\right\|_{\dlp{\infty}}=0.\ee
		In particular, when $\wh{\eta}(0)\ne 0$, we have $\eta=\beta_{\vgu_a,u,\mu}\phi$ and 
		\be\label{conv:vsd:1}\lim_{n\to\infty}\left\|\dm^{|\mu| n}[\sd_{a,\dm I_d}^nv]*u-\beta_{\vgu_a,u,\alpha}\partial^\mu[v*\phi](\dm^{-n}\cdot)\right\|_{\dlp{\infty}}=0.\ee

	\end{enumerate}
\end{theorem}

\bp \begin{enumerate}
	\item[(1)] For every $n\in\N$, define $a_n:=\dm^{-dn}\sd_{a,\dm I_d}^n(\td I_r)$ and $u_n:=a_n*u$. Then 
	$$\wh{a_n}(\xi)=\wh{a}(\dm^{n-1}\xi)\wh{a}(\dm^{n-2}\xi)\dots\wh{a}(\xi),\quad \wh{u_n}(\xi)=\wh{a_n}(\xi)\wh{u}(\xi),\quad\xi\in\dR.$$
	Since $a$ and $u$ are compactly supported, there exists $K\in\N$ such that $\fsupp(a),\fsupp(u)\subseteq[-K,K]^d$, and thus 
	$$\fsupp(a_n), \fsupp(u_n)\subseteq[\dm^{n-1}+\dots+\dm+2][-K,K]^d\subseteq \dm^n[-K,K]^d.$$
	For simplicity, write $\eta^{(\mu)}:=\partial^\mu\eta$ where $\eta$ is the same as in \er{conv:vsd:0}. Then $\eta^{(\mu)}$ must be continuous and supported on $[-K,K]^d$, which further implies that $\wh{\eta^{(\mu)}}$ is an analytic function. Using the Riemann sum, we have
	$$\partial^\nu \wh{\eta^{(\mu)}}(0)=\int_{\dR}\eta^{(\mu)}(x)(-ix)^\nu dx=\lim_{n\to\infty}\dm^{-dn}\sum_{k\in(\dm^{n}[-K,K]^d)\cap\dZ}\eta^{(\mu)}(\dm^{-n}k)(-i\dm^{-n}k)^\nu,$$
	for all $\nu\in\N_0^d$. Define
	$$J_{n,\nu}:=\left[\dm^{-dn}\sum_{k\in(\dm^{n}[-K,K]^d)\cap\dZ}\eta^{(\mu)}(\dm^{-n}k)(-i\dm^{-n}k)^\nu\right]-\partial^\nu \wh{\eta^{(\mu)}}(0).$$
	For any vector $z=(z_1,\dots,z_r)^\tp$ of complex numbers, define $|z|:=|z_1|+\dots+|z_r|$. For $l\in\{0,\dots,m\}$ and $\nu\in\N_0^d$, we have
	\begin{align*}
	&\left|\dm^{(l-|\nu|)n}\partial^\nu\wh{u_n}(0)-\dm^{(l-\tau)n}\partial^\nu\wh{\eta^{(\mu)}}(0)\right|\\
	\le&\left|\dm^{(l-|\nu|)n}\sum_{k\in(\dm^{n}[-K,K]^d)\cap\dZ}u_n(k)(-ik)^\nu-\dm^{(l-\tau-d-|\nu|)n}\sum_{k\in(\dm^{n}[-K,K]^d)\cap\dZ}\eta^{(\mu)}(\dm^{-n}k)(-ik)^\nu\right|+\dm^{(l-\tau)n}|J_{n,\nu}|\\
	\le &\left|\dm^{(l-\tau-d-|\nu|)n}\sum_{k\in(\dm^{n}[-K,K]^d)\cap\dZ}\left[\dm^{\tau n}[\sd_{a,\dm I_d}^n(\td I_r)*u](k)-\eta^{(\mu)}(\dm^{-n}k)\right](-ik)^\nu\right|+\dm^{(l-\tau)n}|J_{n,\nu}|\\
	\le &\left(\dm^{(l-\tau-d-|\nu|)n}\sum_{k\in(\dm^{n}[-K,K]^d)\cap\dZ}|k|^\nu\right)\|\dm^{\tau n}[\sd_{a,\dm I_d}^n(\td I_r)*u]-\eta^{(\mu)}(\dm^{-n}\cdot)\|_{\dlrs{\infty}{r}{1}}+\dm^{(l-\tau)n}|J_{n,\nu}|\\
	\le&\dm^{(l-\tau-d)n}K^{|\nu|d}\|\dm^{\tau n}[\sd_{a,\dm I_d}^n(\td I_r)*u]-\eta^{(\mu)}(\dm^{-n}\cdot)\|_{\dlrs{\infty}{r}{1}}+\dm^{(l-\tau)n}|J_{n,\nu}|
	\end{align*}
	Hence, using \er{conv:vsd:0} and the definition of $J_{n,\nu}$, the above inequalities yield
	$$\lim_{n\to\infty}\left|\dm^{(l-|\nu|)n}\partial^\nu\wh{u_n}(0)-\dm^{(l-\tau)n}\partial^\nu\wh{\eta^{(\mu)}}(0)\right|=0,\quad \forall l\in\N_0\text{ with }l\le \tau,\quad \nu\in\N_0^d,$$
	which further implies
	\be\label{un:eta}\lim_{n\to\infty}\dm^{(l-|\nu|)n}\partial^\nu\wh{u_n}(0)=\lim_{n\to\infty}\dm^{(l-\tau)n}\partial^\nu\wh{\eta^{(\mu)}}(0)=\td(l-\tau)\partial^\nu\wh{\eta^{(\mu)}}(0),\quad \forall l\in\N_0\text{ with }l\le \tau,\quad \nu\in\N_0^d.\ee
	Using \er{sr}, it is not hard to deduce that $\wh{\vgu_a}(\xi)\wh{u}(\xi)=\wh{\vgu_a}(\dm^n\xi)\wh{u_n}(\xi)+\bo(\|\xi\|^{m+1})$ as $\xi\to 0$ for all $n\in\N_0$. Thus for $\nu\in\cup_{q=0}^m\N_{0,q}^d$ such that $|\nu|\le\tau$, using the Leibniz rule and \er{un:eta}, we have
	\begin{align*}
	\partial^\nu[\wh{\vgu_a}\wh{u}](0)&=\lim_{n\to\infty}\partial^\nu[\wh{\vgu_a}(\dm^{n}\cdot)\wh{u_n}](0)=\sum_{\alpha\le\nu}\binom{\nu}{\alpha}\partial^{\nu-\alpha}\wh{\vgu_a}(0)\left[\lim_{n\to\infty}\dm^{(|\nu|-|\alpha|)n}\partial^\alpha\wh{u_n}(0)\right]\\
	&=\lim_{n\to\infty}\dm^{(|\nu|-\tau)n}\sum_{\alpha\le\nu}\binom{\nu}{\alpha}\partial^{\nu-\alpha}\wh{\vgu_a}(0)\partial^\alpha\wh{\eta^{(\mu)}}(0)=\lim_{n\to\infty}\dm^{(|\nu|-\tau)n}\partial^\nu[\wh{\vgu_a}\wh{\eta^{(\mu)}}](0).
	\end{align*}
	Therefore,  \er{vgu:u:0} holds, and this proves item (1).

	\item[(2)]For every $n\in\N_0$, define
	$$J_{n}(\xi):=\wh{\eta^{(\mu)}}(\xi)-\dm^{-dn}\sum_{k\in(\dm^{n}[-K,K]^d)\cap\dZ}\eta^{(\mu)}(\dm^{-n}k)e^{-i(\dm^{-n}k)\cdot\xi},\quad\xi\in\dR.$$
	Using $\fsupp(u_n)\subseteq \dm^n[-K,K]^d$, \er{conv:vsd:0} and the definition of $J_{n}(\xi)$, we deduce that
	\begin{align*}
	\left|\dm^{\tau n}\wh{u_n}(\dm^{-n}\xi)-\wh{\eta^{(\mu)}}(\xi)\right|&\le \dm^{-dn}\left| \sum_{k\in(\dm^{n}[-K,K]^d)\cap\dZ}[\dm^{(\tau+d)n}u_n(k)-\eta^{(\mu)}(\dm^{-n}k)]e^{-i(\dm^{-n}k)\cdot\xi}\right|+|J_n(\xi)|\\
	\le &\dm^{-dn}(2\dm^{n}K+1)^d\|\dm^{(\tau+d)n}u_n-\eta^{(\mu)}(\dm^{-n}\cdot)\|_{\dlrs{\infty}{r}{1}}+|J_n(\xi)|\\
	\le&(3K)^d\|\dm^{\tau n}[\sd_{a,\dm I_r}^n(\td I_r)]-\eta^{(\mu)}(\dm^{-n}\cdot)\|_{\dlrs{\infty}{r}{1}}+|J_n(\xi)|\\
	\to & 0,\quad n\to\infty,\quad\forall\xi\in\dR.
	\end{align*}
	Therefore, \er{ref:eta} must hold. Furthermore,
	\begin{align*}
	\wh{\eta^{(\mu)}}(\dm\xi)&=\lim_{n\to\infty}\left[\dm^{\tau n}\left(\prod_{j=1}^n\wh{a}(\dm^{-(j-1)}\xi)\right)\wh{u}(\dm^{-(n-1)}\xi)\right]\\
	&=\dm^{\tau}\wh{a}(\xi)\lim_{n\to\infty}\left[\dm^{\tau (n-1)}\left(\prod_{j=1}^{n-1}\wh{a}(\dm^{-j}\xi)\right)\wh{u}(\dm^{-(n-1)}\xi)\right]\\
	&=\dm^{\tau}\wh{a}(\xi)\wh{\eta^{(\mu)}}(\xi),\quad\xi\in\dR,
	\end{align*}
	and this proves \er{ref:eta:0}.
	
	\item[(3)]Let $w\in\dlp{0}$. By calculation,
	\begin{align*}
	&\left|\dm^{\tau n}[(\sd_{a,\dm I_d}^n(\td I_r))*(u*w)](k)-\wh{w}(0)\eta^{(\mu)}(\dm^{-n}k)\right|\\
	=&\left|\sum_{z\in\dZ}w(z)\left[\dm^{\tau n}[(\sd_{a,\dm I_d}^n(\td I_r))*u](k-z)-\eta^{(\mu)}(\dm^{-n}k)\right]\right|\\
	\le&\left|\sum_{z\in\dZ}w(z)\left[\dm^{\tau n}[(\sd_{a,\dm I_d}^n(\td I_r))*u](k-z)-\eta^{(\mu)}(\dm^{-n}(k-z))\right]\right|+\left|\sum_{z\in\dZ}w(z)[\eta^{(\mu)}(\dm^{-n}(k-z))-\eta^{(\mu)}(\dm^{-n}k)]\right|\\
	\le&\sum_{z\in\dZ}|w(z)|\left(\|\dm^{\tau n}[\sd_{a,\dm I_d}^n(\td I_r)]*u-\eta^{(\mu)}(\dm^{-n}\cdot)\|_{\dlrs{\infty}{r}{1}}+\sup_{\tilde{k}\in\dZ}|\eta^{(\mu)}(\dm^{-n}\tilde{k}-\dm^{-n}z)-\eta^{(\mu)}(\dm^{-n}\tilde{k})|\right).
	\end{align*}
	Note that all components of $\eta^{(\mu)}$ are continuous and compactly supported, and thus are uniformly continuous, which implies $\lim_{n\to\infty}\sup_{\tilde{k}\in\dZ}|\eta^{(\mu)}(\dm^{-n}\tilde{k}-\dm^{-n}z)-\eta^{(\mu)}(\dm^{-n}\tilde{k})|=0$ for all $z\in\dZ$. Therefore, as $w$ is finitely supported, we conclude from \er{conv:vsd:0} and the above inequalities that item (3) must hold.
\end{enumerate}

Suppose $\tau\in[0,m]$, $\partial^\mu\eta$ is not identically zero and the filter $a$ satisfies \er{eigen:a}. We now prove items (4) and (5).

\begin{enumerate}	
	
	\item[(4)] Since $\eta^{(\mu)}$ is not identically zero and $\wh{\eta^{(\mu)}}(\xi)=(i\xi)^\mu \wh{\eta}(\xi)$, there exists $N\ge m_0$ and $\nu\in\N_{0,N}^d$ with $\nu\ge \mu$ such that $\partial^\alpha\wh{\eta^{(\mu)}}(0)=0$ for all $\alpha\in\cup_{q=0}^{N-1}\N_{0,q}^d$ and $\partial^\nu\wh{\eta^{(\mu)}}(0)\ne 0$. Using \er{ref:eta:0}, we have
	$$\dm^{N}\partial^\nu\wh{\eta^{(\mu)}}(0)=\dm^\tau\sum_{\alpha\le\nu}\binom{\nu}{\alpha}\partial^{\nu-\alpha}\wh{a}(0)\partial^\alpha\wh{\eta^{(\mu)}}(0)=\dm^\tau\wh{a}(0)\partial^\nu\wh{\eta^{(\mu)}}(0),$$
	which implies \be\label{eigen:eta}\dm^{N-\tau}\wh{a}(0)\partial^\nu\wh{\eta^{(\mu)}}(0)=\partial^\nu\wh{\eta^{(\mu)}}(0).\ee
	Therefore, $\dm^{N-\tau}$ is an eigenvalue of $\wh{a}(0)$ and $\partial^\nu\wh{\eta^{(\mu)}}(0)$ is a corresponding eigenvector. Since $N\ge 0$ and $0\le\tau\le m$, by the condition \er{eigen:a}, we must have $\tau=N$ and this proves \er{tau}.\\

	If $\wh{\eta}(0)\ne 0$, let $g_\mu(\xi):=(i\xi)^\mu$, we have
	\be\label{p:mu:al}\partial^\nu\wh{\eta^{(\mu)}}(0)=\sum_{\alpha\le\nu}\binom{\nu}{\alpha}\partial^{\nu-\alpha}g_\mu(0)\partial^\nu\wh{\eta}(0)=\td(\nu-\mu)i^{|\mu|}\mu!\wh{\eta}(0),\quad\forall \nu\in\cup_{q=0}^{|\mu|}\N_{0,q}^d.\ee
	Hence, we conclude that
	\be\label{pd:eta:1}\partial^\mu\wh{\eta^{(\mu)}}(0)=i^{|\mu|}\mu!\wh{\eta}(0)\ne0,\text{ and } \partial^\nu\wh{\eta^{(\mu)}}(0)=0,\quad\forall \nu\in\cup_{q=0}^{|\mu|}\N_{0,q}^d\text{ with }\nu\ne\mu,\ee
	which implies $\tau=|\mu|=m_0$. In this case, we have
	\be\label{pd:eta}\wh{\eta^{(\mu)}}(\xi)=\frac{\partial^\mu\wh{\eta^{(\mu)}}(0)}{\mu!}\xi^\mu+\bo(\|\xi\|^{m_0+1}),\quad\xi\to 0.\ee
	Using \er{ref:eta:0}, $$\dm^{|\mu|}(i\xi)^\mu\wh{\eta}(\dm\xi)=\wh{\eta^{(\mu)}}(\dm\xi)=\dm^\tau\wh{a}(\xi)\wh{\eta^{(\mu)}}(\xi)=\dm^{|\mu|}(i\xi)^\mu \wh{a}(\xi)\wh{\eta}(\xi),\quad\xi\in\dR.$$
	As $\wh{\eta}(0)\ne 0$, so $\eta$ is not identically zero and we must have $\wh{\eta}(\dm\xi)=\wh{a}(\xi)\wh{\eta}(\xi)$. By the eigenvalue condition \er{eigen:a} and the uniqueness of the refinable vector function $\phi$, there must exist $c\in\C\setminus\{0\}$ such that $\eta=c\phi$. This proves item (4).

	\item[(5)] If $\tau>|\mu|=m_0$, then it follows from \er{vgu:u:0} that $\partial^\nu[\wh{\vgu_a}\wh{u}](0)=0$ for all $\nu\in\cup_{q=0}^{m_0}\N_{0,d}^d$. Thus $\beta_{\vgu_a,u,\mu}=0$ and the moment condition \er{mom:vgu:u} clearly holds. If $\tau=|\mu|=m_0$, since $\wh{\eta^{(\mu)}}(\xi)=(i\xi)^\mu\wh{\eta}(\xi)$, we must have \er{pd:eta:1} and \er{pd:eta}. Therefore, we conclude from \er{vgu:u:0}, \er{pd:eta:1} and \er{pd:eta} that \er{mom:vgu:u} must hold.

	For every $v\in\dlrs{0}{1}{r}$, define $\eta_v:=v*\eta:=\sum_{z\in\dZ}v(z)\eta(\cdot-z)$. Then
	$$\partial^\mu\eta_v(\dm^{-n}k)=\sum_{z\in\dZ}v(z)\partial^\mu\eta(\dm^{-n}k-z),\quad\forall k\in\dZ.$$
	By the definition of $\sd_{a,\dm I_d}$, it is easy to verify that
	$$\wh{[(\sd_{a,\dm I_d}^nv)*u]}(\xi)=\dm^{dn}\wh{v}(\dm^n\xi)\wh{a}(\dm^{n-1}\xi)\dots\wh{a}(\xi)\wh{u}(\xi)=\dm^{dn}\wh{v}(\dm^n\xi)\wh{u_n}(\xi)=\wh{\sd_{u_n,\dm^n I_d}v}(\xi),$$
	for all $\xi\in\dR$ and $n\in\N_0$, where $u_n:=\dm^{-dn}[\sd_{a,\dm I_d}^n(\td I_r)]*u$ for all $n\in\N_0$. Hence,
	\begin{align*}
	[(\sd_{a,\dm I_d}^nv)*u](k)=&[\sd_{u_n,\dm^n I_d}v](k)=\sum_{z\in\dZ}v(k)[(\sd_{a,\dm I_d}^n(\td I_r))*u](k-\dm^n z),
	\end{align*}
	for all $k\in\dZ$ and $n\in\N_0$. Therefore, for all $k\in\dZ$, 
	\begin{align*}
	&\left|\dm^{\tau n}[(\sd_{a,\dm I_d}^n v)*u](k)-\partial^\mu\eta_v(\dm^{-n}k)\right|\\
	=&\left|\sum_{z\in\dZ}v(z)\left(\dm^{\tau n}[(\sd_{a,\dm I_d}^n(\td I_r))*u](k-\dm^n z)-\partial^\mu\eta(\dm^{-n}k-z)\right)\right|\\
	\le &\|v\|_{\dlrs{1}{1}{r}}\|\dm^{\tau n}[\sd_{a,\dm I_d}^n(\td I_r)]*u-\partial^\mu\eta(\dm^{-n}\cdot)\|_{\dlrs{\infty}{r}{1}}\\
	\to & 0,\quad n\to\infty,
	\end{align*}
	and this proves \er{conv:vsd:01}.

	Finally, when $\wh{\eta}(0)\ne 0$, then item (4) yields $\tau=|\mu|$,  $\eta=\beta_{\vgu_a,u,\mu}\phi$, and thus  \er{conv:vsd:1} holds.
\end{enumerate}\ep

\section{Convergence of Vector Subdivision Schemes}\label{sec:main}

In this section, we discuss the convergence of a vector subdivision scheme as in Definition~\ref{def:vsd} and some related topics. We will first prove the main result Theorem~\ref{thm:conv:vsd} on the $C^m$-convergence of vector subdivision schemes. Next, we will prove in Theorem~\ref{thm:cr} the convergence rate of a vector subdivision scheme. Then, we will discuss the transformations of vector subdivision schemes using strongly invertible filters to demonstrate the intrinsic connections between different subdivision schemes. Finally, we will apply Theorem~\ref{thm:conv:vsd} to characterize the convergence of Lagrange and Hermite subdivision schemes.

\subsection{Proof of Theorem~\ref{thm:conv:vsd}}

Based on the results of the convergence of vector cascade algorithms from Section~\ref{sec:disc} and with the help of Theorem~\ref{thm:def:vsd}, we are ready to prove Theorem~\ref{thm:conv:vsd} on the convergence of vector subdivision schemes.

\bp[\textbf{Proof of Theorem~\ref{thm:conv:vsd}}]We will prove the theorem by first proving (1) $\Rightarrow$ (2) $\Rightarrow$ (3) $\Rightarrow$ (4), then (4) $\Rightarrow$ (1), and finally (2) $\Leftrightarrow$ (5).\\

\vsp{0.2} 

(1) $\Rightarrow$ (2): As $\sm_\infty(a,\dm I_d)>m$, it follows from Theorem~\ref{thm:vcs} that $\phi\in(C^m(\dR))^{r\times 1}$. As $\phi$ is a fixed point of the vector $\dm I_d$-cascade operator $R_{a,\dm I_d}$, we have 
$$\phi(x)=[R_{a,\dm I_d}^n\phi](x)=\dm^{dn}\sum_{z\in\dZ}a_n(z)\phi(\dm^nx-z),\quad\forall x\in\dR,\quad n\in\N_0,$$
where $a_n:=\dm^{-dn}\sd_{a,\dm I_d}^n(\td I_r)$. Thus for  $\mu\in\cup_{q=0}^{m}\N_{0,q}^d$ and $k\in\dZ$, direct calculation yields
$$\partial^\mu\phi(\dm^{-n}k)=\dm^{(d+|\mu|)n}\sum_{z\in\dZ}a_n(z)\partial^\mu\phi(k-z)=\dm^{(d+|\mu|)n}[a_n*w_\mu](k),$$
where $w_\mu\in\dlrs{0}{r}{1}$ is defined by
\be\label{w:al}w_\mu(k):=\partial^\mu\phi(k),\quad\forall k\in\dZ.\ee
Hence for every $u\in\mom_{\vgu_a,\mu}$, $k\in\dZ$ and $n\in\N_0$, direct calculation yields
\be\label{sd:u}\dm^{|\mu|n}[(\sd_{a,\dm I_d}^n(\td I_r))*u](k)-\beta_{\vgu_a,u,\mu}\partial^\mu\phi(\dm^{-n}k)=\dm^{|\mu|n}\left[(\sd_{a,\dm I_d}^n(\td I_r))*(u-\beta_{\vgu_a,u,\mu}w_\mu)\right](k).\ee
Next, we claim that 
\be\label{u:w:mu}u-\beta_{\vgu_a,u,\mu}w_\mu\in V_{|\mu|,\vgu_a},\quad\forall u\in\mom_{\vgu_a,\mu},\ee
where $V_{|\mu|,\vgu_a}$ is defined as \er{vmy} with $m=|\mu|$ and $y=\vgu_a$. We will prove the claim by showing that
\be\label{vgu:w}\wh{\vgu_a}(\xi)\wh{w_\mu}(\xi)=(i\xi)^\mu+\bo(\|\xi\|^{m+1}),\quad\xi\to0,\quad \forall \mu\in\cup_{q=0}^{m}\N_{0,q}^d,\ee
from which the claim follows trivially. Define $\varphi:=v_a*\phi=\sum_{z\in\dZ}v_a(z)\phi(\cdot-z)$ and 
$$p_\mu(x):=\frac{x^\mu}{\mu!},\quad\forall x\in\dR,\quad \mu\in\cup_{q=0}^m\N_{0,q}^d.$$
Using the Taylor expansion of $p_\mu$, we have
\begin{align*}[p_\mu*\varphi](x):=&\sum_{z\in\dZ}p_\mu(z)\varphi(x-z)=\sum_{\nu\le\mu}\frac{(-1)^{|\nu|}}{\nu!}\partial^\nu p_\mu(x)\sum_{z\in\dZ}(x-z)^\nu\varphi(x-z).
\end{align*}
Note that the function $ \sum_{z\in\dZ}(x-z)^\nu\varphi(x-z)$ is a $\dZ$-periodic function and its Fourier coefficients are given by
\begin{align*}
\int_{[0,1]^d}\sum_{z\in\dZ}(x-z)^\nu\varphi(x-z)e^{-i2\pi k\cdot x}dx=\int_{\dR}x^\nu\varphi(x)e^{-i2\pi k\cdot x}dx=i^{|\nu|}\partial^\nu\wh{\varphi}(2\pi k),
\end{align*}
for all $k\in\dZ$. By item (2) of Proposition~\ref{prop:in}, $\phi$ must satisfy the admissibility condition \er{phi:adm}, which means
$$\partial^\nu\wh{\varphi}(2\pi k)=\td(\nu)\td(k),\quad\forall \nu\in\cup_{q=0}^m\N_{0,q}^d,\quad k\in\dZ,$$
and thus
$$\sum_{z\in\dZ}(x-z)^\nu\varphi(x-z)=i^{|\nu|}\partial^\nu\wh{\varphi}(0)=\td(\nu),\quad \forall \nu\in\cup_{q=0}^m\N_{0,q}^d.$$
Hence, we have
$$[p_\mu*\varphi](x)=\sum_{\nu\le\mu}\frac{(-1)^{|\nu|}}{\nu!}\partial^\nu p_\mu(x)\td(\nu)=p_\mu(x),\quad \forall x\in\dR,\quad\mu\in\cup_{q=0}^m\N_{0,q}^d.$$
Now, for every $\nu,\mu\in\N_0^d$ and $k\in\dZ$, we have
\begin{align*}
\partial^\mu p_{\nu}(k)=&\sum_{\alpha\le\nu}p_\nu(k-z)\partial^\mu\varphi(z)=\sum_{\alpha\le\nu}\frac{(-1)^{|\alpha|}}{\alpha!}\partial^\alpha p_\nu(k)\sum_{z\in\dZ}z^\alpha\partial^\mu\varphi(z)\\
=&\sum_{\alpha\le\nu}\frac{(-1)^{|\alpha|}}{\alpha!}\partial^\alpha p_\nu(k)\sum_{z\in\dZ}z^\alpha[\vgu_a*\partial^\mu\phi](z)=\sum_{\alpha\le\nu}\frac{(-1)^{|\alpha|}}{\alpha!}\partial^\alpha p_\nu(k)\sum_{z\in\dZ}z^\alpha[\vgu_a*w_\mu](z)\\
=&\sum_{\alpha\le\nu}\frac{(-i)^{|\alpha|}}{\alpha!}\partial^\alpha p_\nu(k)\partial^\alpha[\wh{(\vgu_a*w_\mu)}](0)=\sum_{\alpha\le\nu}\frac{(-i)^{|\alpha|}}{\alpha!}\partial^\alpha p_\nu(k)\partial^\alpha[\wh{\vgu_a}\wh{w_\mu}](0)
\end{align*}
From the above identities, for $\nu,\mu\in\cup_{q=0}^m\N_{0,q}^d$:
\begin{itemize}
	\item If $\mu\le\nu$, as all $\partial^\alpha p_\nu$, $\alpha\le\nu$ are linearly independent, then the above identities yield
	$$\partial^\alpha[\wh{\vgu_a}\wh{w_\mu}](0)=\td(\alpha-\mu)i^{|\mu|}\mu!,\quad\forall \alpha\in\N_0^d\text{ with }\alpha\le\nu;$$
	
	\item Otherwise, if $\mu\le\nu$ does not hold, then there must exist $j\in\{1,\dots,d\}$ such that $\mu_j>\nu_j$, where $\mu_j$ and $\nu_j$ are the $j$-th components of $\mu$ and $\nu$ respectively. Hence, we must have $\partial^\mu p_\nu\equiv 0$ and the linear independence of $\partial^\alpha p_\nu$, $\alpha\le\nu$ yields
	$$\partial^\alpha[\wh{\vgu_a}\wh{w_\mu}](0)=0,\quad\forall \alpha\in\N_0^d\text{ with }\alpha\le\nu;$$
\end{itemize}
Therefore, we have
$$\partial^\nu[\wh{\vgu_a}\wh{w_\mu}](0)=\td(\nu-\mu)i^{|\mu|}\mu!,\quad\forall \nu, \mu\in\cup_{q=0}^{m}\N_{0,q}^d,$$
which is equivalent to \er{vgu:w}, and this proves the claim.

Finally, since $u-\beta_{\vgu_a,u,\mu}w_\mu\in V_{|\mu|,\vgu_a}$ for all $u\in\mom_{\vgu_a,\mu}$, item (2) follows immediately from \er{sd:u} and item (3) of Theorem~\ref{thm:vcs}.\\

\vsp{0.2}

(2) $\Rightarrow$ (3) $\Rightarrow$ (4): Trivial.\\

\vsp{0.2}
(4) $\Rightarrow$ (1): Suppose $u\in\mom_{\vgu_a,\mu}$, then there exist $u_1,\dots,u_t\in G_{\vgu_a,\mu}$ and $c_1,\dots,c_t\in\dlp{0}$ such that $u=u_1*c_1+\dots+u_t*c_t$. Since $u_1,\dots,u_t\in\mom_{\vgu_a,\mu}$,  for any $\nu\in\N_0^d$ with $|\nu|\le|\mu|$, we have
\begin{align*}\partial^\nu[\wh{\vgu_a}\wh{u_j}\wh{c_j}](0)=&\sum_{\alpha\le\nu}\binom{\nu}{\alpha}\partial^{\nu-\alpha}[\wh{\vgu_a}\wh{u_j}](0)\partial^{\alpha}\wh{c_j}(0)=\td(\mu-\nu)\partial^\mu[\wh{\vgu_a}\wh{u_j}](0)\wh{c_j}(0),\quad\forall j=1,\dots,t.
\end{align*}
Thus $u_1*c_1,\dots,u_t*c_t\in \mom_{\vgu_a,\mu}$ and
\begin{align*}
\beta_{\vgu_a,u,\mu}=&\frac{\partial^{\mu}[\wh{\vgu_a}\wh{u}](0)}{i^{|\mu|}\mu!}=\frac{\partial^\mu[\wh{\vgu_a}\wh{u_1}\wh{c_1}](0)}{i^{|\mu|}\mu!}+\dots+\frac{\partial^\mu[\wh{\vgu_a}\wh{u_t}\wh{c_t}](0)}{i^{|\mu|}\mu!}\\
=&\frac{\partial^\mu[\wh{\vgu_a}\wh{u_1}](0)\wh{c_1}(0)}{i^{|\mu|}\mu!}+\dots+\frac{\partial^\mu[\wh{\vgu_a}\wh{u_t}](0)\wh{c_t}(0)}{i^{|\mu|}\mu!}=\beta_{\vgu_a,u_1,\mu}\wh{c_1}(0)+\dots+\beta_{\vgu_a,u_t,\mu}\wh{c_t}(0).
\end{align*}
Hence
\begin{align*}
&\|\dm^{|\mu|n}[\sd_{a,\dm I_d}^n(\td I_r)]*u-\beta_{\vgu_a,u,\mu}\partial^\mu\phi(\dm^{-n}\cdot)\|_{\dlrs{\infty}{r}{1}}\\
=&\|\dm^{|\mu|n}[\sd_{a,\dm I_d}^n(\td I_r)]*(u_1*c_1+\dots+u_t*c_t)-[\beta_{\vgu_a,u_1,\mu}\wh{c_1}(0)+\dots+\beta_{\vgu_a,u_t,\mu}\wh{c_t}(0)]\partial^\mu\phi(\dm^{-n}\cdot)\|_{\dlrs{\infty}{r}{1}}\\
\le &\sum_{j=1}^t\|\dm^{|\mu|n}[\sd_{a,\dm I_d}^n(\td I_r)]*(u_j*c_j)-\beta_{\vgu_a,u_j,\mu}\wh{c_j}(0)\partial^\mu\phi(\dm^{-n}\cdot)\|_{\dlrs{\infty}{r}{1}}.
\end{align*}
Since $u_j\in G_{\vgu_a,\mu}$ for all $j=1,\dots,t$, we have
$$\lim_{n\to\infty}\sum_{j=1}^t\|\dm^{|\mu|n}[\sd_{a,\dm I_d}^n(\td I_r)]*u_j-\beta_{\vgu_a,u_j,\mu}\partial^\mu\phi(\dm^{-n}\cdot)\|_{\dlrs{\infty}{r}{1}}=0,$$
and thus item (3) of Theorem~\ref{thm:def:vsd} yields 
$$\lim_{n\to\infty}\|\dm^{|\mu|n}[\sd_{a,\dm I_d}^n(\td I_r)]*(u_j*c_j)-\beta_{\vgu_a,u_j,\mu}\wh{c_j}(0)\partial^\mu\phi(\dm^{-n}\cdot)\|_{\dlrs{\infty}{r}{1}}=0,$$
for all $j=1,\dots,t$. Therefore, \er{conv:vsd:phi} holds for all $u\in \mom_{\vgu_a,\mu}$.\\

Define $V_{m,\vgu_a}$ via \er{vmy} with $y=\vgu_a$. Then it is clear that $V_{m,\vgu_a}\subseteq \mom_{\vgu_a,\mu}$ and
$$\beta_{\vgu_a,u,\mu}=\partial^\mu[\wh{\vgu_a}\wh{u}](0)=0,\quad\forall u\in V_{m,\vgu_a}.$$
Now, using the above identity, \er{conv:vsd:phi} and $|\mu|=m$, we have
$$\lim_{n\to\infty}\|\dm^{mn}[\sd_{a,\dm I_d}^n(\td I_r)]*u\|_{\dlrs{\infty}{r}{1}}=0,\quad\forall u\in V_{m,\vgu_a},$$
that is, item (3) of Theorem~\ref{thm:vcs} holds with $p=\infty$. Therefore, by Theorem~\ref{thm:vcs}, item (1) holds.\\

\vsp{0.2}

(2) $\Leftrightarrow$ (5): Suppose item (2) holds. Let $v\in\dlrs{0}{1}{r}$ be an input data. By using an argument similar to the proof of \er{conv:vsd:1} in item (5) of Theorem~\ref{thm:def:vsd}, we can deduce from \er{conv:vsd:phi} that \er{conv:vsd} must hold with $\eta_v:=v*\phi$ and this proves the implication (2) $\Rightarrow$ (5).\\

Conversely, suppose item (5) holds, it follows from \er{conv:vsd} that there exists $\eta\in (C^m(\dR))^{r\times 1}$ that satisfies
\be\label{sd:eta:0}\lim_{n\to\infty}\|\dm^{|\mu|n}[\sd_{a,\dm I_d}^n(\td I_r)]*u-\beta_{\vgu_a,u,\mu}\partial^\mu \eta(\dm^{-n}\cdot)\|_{\dlrs{\infty}{r}{1}}=0,\quad \forall u\in\mom_{\vgu_a,\mu},\quad  \mu\in\cup_{q=0}^m\N_{0,q}^d.\ee
Then by item (1) of Theorem~\ref{thm:def:vsd}, \er{sd:eta:0} yields
\begin{align*}\partial^\mu[\wh{\vgu_a}\wh{u}](0)=&\beta_{\vgu_a,u,\mu}\partial^\mu[\wh{\vgu_a}\wh{\partial^\mu\eta}](0)=\beta_{\vgu_a,u,\mu}\sum_{\nu\le\mu}\binom{\mu}{\nu}\partial^{\mu-\nu}\wh{\vgu_a}(0)\partial^{\nu}\wh{\partial^\mu\eta}(0)=\beta_{\vgu_a,u,\mu}\wh{\vgu_a}(0)\partial^\mu\wh{\partial^\mu\eta}(0)\\
=&\beta_{\vgu_a,u,\mu}i^{|\mu|}\mu!\wh{\vgu_a}(0)\wh{\eta}(0)=\partial^\mu[\wh{\vgu_a}\wh{u}](0)\wh{\vgu_a}(0)\wh{\eta}(0),\quad\forall u\in\mom_{\vgu_a,\mu}.\end{align*}
Hence, $\wh{\vgu_a}(0)\wh{\eta}(0)=1$, which, in particular, implies that $\wh{\eta}(0)\ne 0$. Therefore, by item (5) of Theorem~\ref{thm:def:vsd}, we must have $\eta=\phi$, and this proves the implication (5) $\Rightarrow$ (2). \\

\vsp{0.2}

Finally, if $\phi\in (C^m(\dR))^{r\times 1}$ has stable integer shifts, we conclude from Theorem~\ref{thm:vcs} that all items (1)-(5) hold. \ep

\subsection{Convergence Rate of a Vector Subdivision Scheme}

Using a vector subdivision scheme (more specifically, item (2) of Theorem~\ref{thm:conv:vsd}), we have a method to compute the refinable function vector $\phi$ and its derivatives. Our next task is to analyze how fast a vector subdivision scheme converges, essential for computation efficiency in practice. This leads to the convergence rate of a vector subdivision scheme, and we have the following result.

\begin{theorem}\label{thm:cr}Let $\dm\in\setminus\{1\}$, $r\in\N$ and $m\in\N_0$. Let $a\in\dlrs{0}{r}{r}$ be such that $m<\sm_\infty(a,\dm I_d)\le m+1$ and $\vgu_a\in\dlrs{0}{1}{r}$ is an order $m+1$ matching filter of $a$ that satisfies \er{sr}. Then $a$ must satisfy \er{eigen:a}, which means the set $\{\partial^\mu\vgu_a(0):\,\mu\in\cup_{q=0}^m\N_{0,q}^d\}$ is uniquely determined through \er{sr:diff:0} and the standard $\dm I_d$-refinable function vector $\phi$ of $a$ must satisfy $\phi\in (C^m(\dR))^{r\times 1}$. Let $\mu\in\cup_{q=0}^m\N_{0,q}^d$ and $u\in\mom_{\vgu_a,\mu}$. Define
	\be\label{drv}\begin{aligned}\drv(\vgu_a*u):=\inf\{&N\in\N:\,\partial^\nu[\wh{\vgu_a}\wh{u}-\beta_{\vgu_a,u,\mu}(i\cdot)^\mu](0)=0\text{ for all }\nu\in\cup_{q=0}^{N-1}\N_{0,q}^d\\
		&\text{and }\partial^\nu[\wh{\vgu_a}\wh{u}-\beta_{\vgu_a,u,\mu}(i\cdot)^\mu](0)\ne0\text{ for some }\nu\in \N_{0,N}^d\}>|\mu|,\end{aligned}\ee
	and $S_{\mu,u}:=\min\{\drv(\vgu_a*u)-|\mu|,\,\sm_\infty(a,\dm I_d)-|\mu|\}>0$. If $\drv(\vgu_a*u)\le m+1$, then for every $\veps\in(0,S_{\mu,u})$, there exists a constant $C>0$ such that
	\be\label{cr:vsd}\|\dm^{|\mu|n}[\sd_{a,\dm I_d}^n(\td I_r)]*u-\beta_{\vgu_a,u,\mu}\partial^\mu\phi(\dm^{-n}\cdot)\|_{\dlrs{\infty}{r}{1}}\le C\dm^{-(S_{\mu,u}-\veps)n},\quad\forall n\in\N,\ee
	and for every $v\in\dlrs{0}{1}{r}$, we have
	\be\label{cr:vsd:1}\|\dm^{|\mu|n}[\sd_{a,\dm I_d}^nv]*u-\beta_{\vgu_a,u,\mu}\partial^\mu[v*\phi](\dm^{-n}\cdot)\|_{\dlp{\infty}}\le C_v\dm^{-(S_{\mu,u}-\veps)n},\quad\forall n\in\N,\ee
	where $C_v:=\|v\|_{\dlrs{1}{1}{r}}$.
	
\end{theorem}

\bp First note that by Theorem~\ref{thm:vcs}, we have $m_a:=\sr(a,\dm I_d)\ge\sm_\infty(a,\dm I_d)>m$ and $a$ must satisfy \er{eigen:a}. Thus there must exist an order $m+1$ matching filter $\vgu_a\in\dlrs{0}{1}{r}$ that satisfies \er{sr} and a unique vector $\phi$ of compactly supported distributions that satisfies \er{stan:vref}. Furthermore, $\sm_\infty(a,\dm I_d)$ yields $\phi\in (C^m(\dR))^{r\times1}$ and Proposition~\ref{prop:in} tells us that $\phi$ must satisfy the admissibility condition \er{phi:adm}. For each $\mu\in \cup_{q=0}^m\N_{0,q}^{d}$ and define $w_\mu\in\dlrs{0}{r}{1}$ via \er{w:al}. By using the same argument as in the proof of (1) $\Rightarrow$ (2) in Theorem~\ref{thm:conv:vsd}, one can show that \er{u:w:mu} must hold, which further implies
$$\wh{\vgu_a}(\xi)[\wh{u}(\xi)-\beta_{\vgu_a,u,\mu}\wh{w_{\mu}}](\xi)=\bo(\|\xi\|^{\drv(\vgu_a*u)}),\quad\xi\to 0$$
for all $u\in\mom_{\vgu_a,\mu}$. Hence for every $u\in\mom_{\vgu_a,\mu}$, we have \be\label{u:w:mu:0}u-\beta_{\vgu_a,u,\mu}w_{\mu}\in V_{\drv(\vgu_a*u)-1,\vgu_a}.\ee
For every $\veps\in(0,S_{\mu,u})$, by the definition of $\rho_{\drv(\vgu_a*u)}(a,\dm I_d,\vgu_a)_\infty$, there exists a constant $C>0$ such that
\be\label{rho:s}\|\dm^{|\mu|n}[\sd_{a,\dm I_d}^n(\td I_r)]*(u-\beta_{\vgu_a,u,\mu}w_\mu)\|_{\dlrs{\infty}{r}{1}}\le C\dm^{|\mu|n}\dm^{\veps n}[\rho_{\drv(\vgu_a*u)}(a,\dm I_d,\vgu_a)_\infty]^n,\quad\forall n\in\N.\ee
Now suppose $\drv(\vgu_a*u)\le m+1$. If $\drv(\vgu_a*u)=m+1$, then using $\sm_{\infty}(a,\dm I_d)\le m+1\le m_a$ and Theorem~\ref{thm:sm}, we have
$$\rho_{\drv(\vgu_a*u)}(a,\dm I_d,\vgu_a)_\infty=\rho_{m+1}(a,\dm I_d,\vgu_a)_\infty=\dm^{-\sm_\infty(a,\dm I_d)}.$$
If $\drv(\vgu_a*u)< m+1$, then $\drv(\vgu_a*u)\le m$ as $\drv(\vgu_a*u)$ is an integer. Using Theorem~\ref{thm:sm} and $\sm_\infty(a,\dm I_d)>m$, we have
$$\rho_{m_a}(a,\dm I_d,\vgu_a)_\infty=\dm^{-\sm_\infty(a,\dm I_d)}<\dm^{-m}\le\dm^{-\drv(\vgu_a*u)},$$
and it then follows from \er{rho:ko:ma} that
$$\rho_{\drv(\vgu_a*u)}(a,\dm I_d,\vgu_a)_\infty=\max\{\rho_{m_a}(a,\dm I_d,\vgu_a)_\infty,\dm^{-\drv(\vgu_a*u)}\}=\dm^{-\drv(\vgu_a*u)}.$$
Therefore, we always have
\be\label{rho:s:0}\rho_{\drv(\vgu_a*u)}(a,\dm I_d,\vgu_a)_\infty\le\dm^{-\min\{\sm_\infty(a,\dm I_d),\,\drv(\vgu_a*u)\}}=\dm^{-S_{\mu,u}-|\mu|}.\ee
Consequently, we conclude from \er{rho:s} and \er{rho:s:0} that
\begin{align*}
&\|\dm^{|\mu|n}[\sd_{a,\dm I_d}^n(\td I_r)]*(u-\beta_{\vgu_a,u,\mu}w_\mu)\|_{\dlrs{\infty}{r}{1}}\le C\dm^{-(S_{\mu,u}-\veps)n},\quad\forall n\in\N.
\end{align*}
Now \er{cr:vsd} follows immediately from the above inequalities and \er{sd:u}. Finally, for every $v\in\dlrs{0}{1}{r}$, by using an argument similar to the proof of \er{conv:vsd:1} in item (5) of Theorem~\ref{thm:def:vsd}, one can conclude that 
$$\left\|\dm^{|\mu|n}[\sd_{a,\dm I_d}^n v]*u-\partial^\mu[v*\phi](\dm^{-n}\cdot)\right\|_{\dlp{\infty}}\le \|v\|_{\dlrs{1}{1}{r}}\|\dm^{|\mu|n}[\sd_{a,\dm I_d}^n(\td I_r)]*u-\partial^\mu\phi(\dm^{-n}\cdot)\|_{\dlrs{\infty}{r}{1}},$$
for all $n\in\N$ and thus \er{cr:vsd:1} follows trivially.  \ep

\subsection{Transformations of Vector Subdivision Schemes Using Strongly Invertible Filters}

In this subsection, we show that every vector subdivision scheme can be transformed into a scalar-like subdivision scheme. For the scalar case $r=1$, we see from Definition~\ref{def:svd:0} that a matching filter is not involved for the convergence of a scalar subdivision scheme, and let us first explain why this is the case. Let $\dm\in\N\setminus\{1\}$ and suppose $a\in\dlp{0}$ with $\wh{a}(0)=1$,  then clearly the condition \er{eigen:a} holds. In this case, it is well-known that there is a unique compactly supported distribution $\phi$ such that $\wh{\phi}(\dm\xi)=\wh{a}(\xi)\phi(\xi)$ and $\wh{\phi}(0)=1$, which is known as the \emph{standard $\dm I_d$-refinable function of the filter $a$} and is defined as
\be\label{stan:ref}\wh{\phi}(\xi):=\prod_{j=1}^\infty\wh{a}(\dm^{-j}\xi),\quad\forall\xi\in\dR.\ee
With $\phi$ defined as above, we see that
$$\frac{1}{\wh{\phi}(\dm\xi)}\wh{a}(\xi)=\frac{1}{\wh{\phi}(\xi)},\quad\forall\xi\in\dR,$$
and any order $m+1$ matching filter $\vgu_a\in\dlp{0}$ in \er{sr} must satisfy
\be\label{vgu:phi}\wh{\vgu_a}(\xi)=\frac{1}{\wh{\phi}(\xi)}+\bo(\|\xi\|^{m+1}),\quad\xi\to 0.\ee
From \er{vgu:phi}, we must have $\wh{\vgu_a}(0)=1$ and thus the moment condition in \er{sr} simplifies to
\be\label{sr:s}\wh{a}(\xi+2\pi\omega)=\bo(\|\xi\|^{m+1}),\quad\xi\to 0,\quad\forall\omega\in\Omega_{\dm I_d}\setminus\{0\},\ee
which does not explicitly involve the matching filter $\vgu_a$. Consequently, we do not need a matching filter for the sum rule condition \er{sr} of a scalar filter $a$. Furthermore, for all $\mu\in\cup_{q=0}^m\N_{0,q}^d$, the set $\mom_{\vgu_a,\mu}$ is given by
$$\mom_{\vgu_a,\mu}=\{u\in\dlp{0}:\,\wh{u}(\xi)=\wh{\phi}(\xi)\beta_{u,\phi,\mu}(i\xi)^\mu+\bo(\|\xi\|^{|\mu|+1}),\,\xi\to 0\},\text{ where }\beta_{u,\phi,\mu}:=\frac{\partial^\mu[\wh{u}/\wh{\phi}](0)}{i^{|\mu|}\mu!},$$
which does not explicitly involve the matching filter $\vgu_a$. Therefore, according to Definition~\ref{def:svd:0}, we do not need a matching filter for the convergence of a scalar subdivision scheme.

Now, consider the case $r>1$. Let $\dm\in\N\setminus\{1\}$, $m\in\N_0$ and $a\in\dlrs{0}{r}{r}$. If $a$ has order $m+1$ sum rules with respect to $\dm I_d$ such that \er{sr} holds with a matching filter $\vgu_a=\td [\pe_1^r]^\tp=[\td,0,\dots,0]\in\dlrs{0}{1}{r}$ (or equivalently, $\wh{\vgu_a}(\xi)=[\pe_1^r]^\tp=[1,0,\dots,0]$), then $a$ is called an \emph{order $m+1$ scalar-type filter}. In this case, let $a=\begin{bmatrix}a_{1,1} & a_{1,2}\\
a_{2,1} & a_{2,2}\end{bmatrix}$ where $a_{1,1}\in\dlp{0}$, $a_{1,2}\in\dlrs{0}{1}{(r-1)}$, $a_{2,1}\in\dlrs{0}{(r-1)}{1}$ and $a_{2,2}\in\dlrs{0}{(r-1)}{(r-1)}$, it is easy to deduce from \er{sr} that $a$ must satisfy
\be\label{nf:a}\wh{a_{1,1}}(\xi+2\pi\omega)=\td(\omega)+\bo(\|\xi\|^{m+1}),\quad \wh{a_{1,2}}(\xi+2\pi\omega)=\bo(\|\xi\|^{m+1}),\quad \xi\to 0\quad\forall\omega\in\Omega_{\dm I_d}.\ee
The moment conditions in \er{nf:a} are very similar to the sum rule condition \er{sr:s} of a scalar filter, and the matrix filter $a$ essentially behaves like a scalar filter $a_{1,1}$. When $a$ is a scalar-type filter, techniques from the scalar case $r=1$ can be borrowed or generalized directly to study the properties of $a$ and the subdivision scheme that uses $a$.

If a vector subdivision scheme uses a scalar-type filter, we call it \emph{a scalar-type vector subdivision scheme}. We claim that every vector subdivision scheme can be transformed into a scalar-type one, and thus, the two subdivision schemes have no fundamental difference. We say that a filter $U\in\dlrs{0}{r}{r}$ is \emph{strongly invertible} if $\det(\wh{U}(\xi))=c_0e^{-ik\cdot\xi}$ for some $c_0\in\C\setminus\{0\}$ and $k\in\Z$ for all $\xi\in\dR$. It is easy to see that $U\in\dlrs{0}{r}{r}$ is strongly invertible if and only if there exists $V\in\dlrs{0}{r}{r}$ such that $\wh{V}=\wh{U}^{-1}$.  Here, we recall the following result regarding transformations with strongly invertible matrices.
\begin{lemma}\label{thm:nf}[\cite[Lemma 2.3]{han10}, also see \cite[Lemma 3.1]{hl20pp}] Suppose $v,w\in\dlrs{0}{s}{r}$ satisfies $\wh{v}(0)\ne 0$ and $\wh{w}(0)\ne 0$. Then for any $n\in\N_0$, there exists a strongly invertible filter $U_n\in\dlrs{0}{r}{r}$ such that
	$$\wh{v}(\xi)=\wh{w}(\xi)\wh{U_n}(\xi)+\bo(\|\xi\|^{n}),\quad\xi\to 0.$$
\end{lemma}
With the help of Lemma~\ref{thm:nf}, we now show how to transform a vector subdivision scheme to a scalar-type one and, more importantly, the intrinsic connections between the original vector subdivision scheme and the one we get after transformation. Let $\dm\in\N\setminus\{1\}$ and $m\in\N_0$. Suppose $a\in\dlrs{0}{r}{r}$ has order $m+1$ sum rules with respect to $\dm I_d$ and $\vgu_a\in\dlrs{0}{1}{r}$ is a matching filter that satisfies \er{sr}. Further, assume that $a$ satisfies \er{eigen:a}. By Lemma~\ref{thm:nf}, there exists a strongly invertible filter $U\in\dlrs{0}{r}{r}$ such that
$$\wh{\vgu_a}(\xi)\wh{U}(\xi)=[1,0,\dots,0]+\bo(\|\xi\|^{m+1}),\quad\xi\to 0.$$
Let $\phi$ be the standard $\dm I_d$-refinable function vector of $a$ that satisfies \er{stan:vref}. Define a new filter $\mra\in\dlrs{0}{r}{r}$ and a new vector $\mrphi$ of compactly supported distributions via
\be\label{mr:a:phi}\wh{\mra}(\xi):=\wh{U}(\dm\xi)^{-1}\wh{a}(\xi)\wh{U}(\xi),\quad\wh{\mrphi}(\xi):=\wh{U}(\xi)^{-1}\wh{\phi}(\xi),\quad\forall\xi\in\dR.\ee
Then we obtain the following:
\begin{enumerate}
	\item[(1)] $\wh{\mrphi}(\dm\xi)=\wh{\mra}(\xi)\wh{\mrphi}(\xi)$ and $[1,0,\dots,0]\wh{\mrphi}(0)=1$.
	
	\item[(2)] $\mra$ has order $m+1$ sum rules with respect to $\dm I_d$ with a matching filter $\vgu_{\mra}:=\td [\pe_1^r]^\tp=[\td,0,\dots,0]$, that is, \er{sr} holds with $a$ being replaced by $\mra$. Therefore, $\mra$ is an order $m+1$ scalar-type filter.
	
	\item[(3)]Recall the definition of $V_{m,y}$ as in \er{vmy} for all $y\in\dlrs{0}{1}{r}$ with $\wh{y}(0)\ne 0$. For every $j\in\{0,\dots,m+1\}$ and $\mu\in\cup_{q=0}^m\N_{0,q}^d$, the following maps are bijective:
	$$f_j:\,V_{j-1,\td [\pe_1^r]^\tp}\to V_{j-1,\vgu_a},\quad f_j(w)=U*w,$$
	$$g_\mu:\,\mom_{\td [\pe_1^r]^\tp,\mu}\to \mom_{\vgu_a,\mu},\quad g_\mu(u)=U*u.$$
	Consequently, the set
	$$B_{j-1,\vgu_a}:=\left\{U*\nabla^\nu\td \pe_1^r:\,\nu\in\cup_{q=0}^j\N_{0,q}^d\right\}\cup\{U*\td \pe_l^r:\,l=2,\dots,r\}$$
	generates $V_{j-1,\vgu_a}$ for every $j\in\{0,\dots,m+1\}$, and the set
	$$G_{\vgu_a,\mu}:=\left\{U*\nabla^\mu\td \pe_1^r\right\}\cup\{U*\td \pe_l^r:\,l=2,\dots,r\}$$
	generates $\mom_{\vgu_a,\mu}$ for every $\mu\in\cup_{q=0}^m\N_{0,q}^d$.
	
	\item[(4)]For every $n\in\N_0$, define $a_n:=\dm^{-dn}\sd_{a,\dm I_d}^n(\td I_r)$ and $\mra_n:=\dm^{-dn}\sd_{\mra,\dm I_d}^n(\td I_r)$. For every $j\in\{0,\dots,m+1\}$, let $\mrw\in V_{j-1,\td [\pe_1^r]^\tp}$ and define $w:=U*\mrw\in V_{j-1,\vgu_a}$, then 
	$$\wh{a_n}(\xi)\wh{w}(\xi)=\wh{U}(\dm^{n}\xi)\wh{\mra_n}(\xi)\wh{\mrw}(\xi),\quad\forall\xi\in\dR,$$
	or equivalently,
	$$[\sd_{a,\dm I_d}^n(\td I_r)]*w=[\sd_{\mra,\dm I_d}^nU]*\mrw=\sum_{z\in\dZ}U(z)[(\sd_{\mra,\dm I_d}^n(\td I_r))*\mrw](\cdot-\dm^nz).$$
	For $p\in[1,\infty]$, we deduce from the above identity that
	\begin{align*}\limsup_{n\to\infty}\|[\sd_{a,\dm I_d}^n(\td I_r)]*w\|_{\dlrs{p}{r}{1}}^{\frac{1}{n}}\le&\limsup_{n\to\infty}\|[\sd_{\mra,\dm I_d}^n(\td I_r)]*\mrw\|_{\dlrs{p}{r}{1}}^{\frac{1}{n}},\quad\forall \mrw\in V_{j-1,\td [\pe_1^r]^\tp},\end{align*}
	which implies $\rho_j(a,\dm I_d,\vgu_a)_p\le\rho_j(\mra,\dm I_d,\td [\pe_1^r]^\tp)_p$ for all $j\in\{0,\dots,m+1\}.$. Similarly, one can show that $\rho_j(\mra,\dm I_d,\td[\pe_1^r]^\tp)_p\le\rho_j(a,\dm I_d,\vgu_a)_p$ for all $j\in\{0,\dots,m+1\}.$ Therefore, we conclude that
	$$\rho(a,\dm I_d)_p=\rho(\mra,\dm I_d)_p,\quad \sm_p(a,\dm I_d)=\sm_p(\mra,\dm I_d),\quad 1\le p\le\infty.$$
	Consequently, using Theorem~\ref{thm:conv:vsd}, the vector subdivision scheme with the filter $a$ is $C^m$-convergent if and only if the scalar-type vector subdivision scheme with the filter $\mra$ is $C^m$-convergent.
	
\end{enumerate}

\subsection{Convergence of Lagrange and Hermite Subdivision Schemes}

We now apply our theory on vector subdivision schemes to discuss Lagrange and Hermite subdivision schemes. Both schemes are the special cases of the generalized Hermite subdivision schemes. Here, we first recall the convergence of a generalized Hermite subdivision scheme recently introduced in \cite{han23-2}.

\begin{definition}\label{def:ghsd}Let $r\in\N$, $m\in\N_0$ and let $\Lambda:=\{\nu^1,\dots,\nu^r\}\subseteq\cup_{q=0}^m\N_{0,q}^d$ be an ordered multiset with $\nu^1=0$. A \emph{generalized Hermite subdivision scheme of type $\Lambda$ with filter $a\in\dlrs{0}{r}{r}$} is \emph{convergent with limit functions in $C^m(\dR)$} if for every input data $v\in\dlrs{0}{1}{r}$, there exists a limit function $\eta_v\in C^m(\dR)$ such that 
	\be\label{conv:ghsd}\lim_{n\to\infty}\|2^{|\nu^l|n}[\sd_{a,2 I_d}^nv]\pe_l^r-\partial^{\nu^l}\eta_v(2^{-n}\cdot)\|_{\dlp{\infty}}=0,\quad\forall l=1,\dots,r.\ee
	In particular:
	\begin{itemize}
		\item A \emph{Lagrange subdivision scheme} is a generalized Hermite subdivision scheme of type $\Lambda:=\{0,0,\dots,0\}$, that is, $\nu^l=0$ for all $l=1,\dots,r$.
		
		\item When $d=1$, a \emph{Hermite subdivision scheme or order $r$} is a generalized Hermite subdivision scheme of type $\Lambda:=\{0,1,\dots,r-1\}$.
	\end{itemize}
	
\end{definition}

Comprehensive structural characterization and convergence analysis of a generalized Hermite subdivisions scheme has been established in \cite[Theorems 1.1 and 1.2]{han23-2}, with the extra assumption that the refinable function $\phi$ of the filter $a$ has stable integer shifts. Using our theory of subdivision schemes from the previous sections, we can drop the stability on the integer shifts of $\phi$  and obtain the following result that improves \cite[Theorems 1.1 and 1.2]{han23-2}.

\begin{theorem}\label{thm:ghsd}Let $r\in\N$, $m\in\N_0$ and let $\Lambda:=\{\nu^1,\dots,\nu^r\}\subseteq\cup_{q=0}^m\N_{0,q}^d$ be an ordered multiset with $\nu^1=0$. Define $\ttm:=\max\{|\nu^1|,\dots,|\nu^r|\}$. Let $a\in\dlrs{0}{r}{r}$ be such that
	\be\label{eigen:a:1}1\text{ is a simple eigenvalue of }\wh{a}(0)\text{ and }2^{n}\text{ is not an eigenvalue of }a\text{ for all }n\in \N.\ee 
	Suppose $a$ has order $\ttm+1$ sum rules with respect to $2 I_d$, $\vgu_a\in\dlrs{0}{1}{r}$ is an order $\ttm+1$ matching filter of $a$, and $\phi$ is the unique $r\times 1$ vector of compactly supported distributions that satisfies \er{stan:vref} with $\dm=2$. Then the following statements hold:
	\begin{enumerate}
		\item[(1)]If the generalized Hermite subdivision scheme of type $\Lambda$ with the filter $a$ is convergent with limit functions in $C^m(\dR)$, then $\sm_\infty(a,2I_d)>0$, $\phi\in (C^m(\dR))^{r\times1}$ and up to multiplying a non-zero constant, $\vgu_a$ must satisfy
		\be\label{ghsd:vgu}\wh{\vgu_a}(\xi)=[(i\xi)^{\nu^1}+\bo(\|\xi\|^{|\nu^1|+1}),\,(i\xi)^{\nu^2}+\bo(\|\xi\|^{|\nu^2|+1}),\,\dots\,(i\xi)^{\nu^r}+\bo(\|\xi\|^{|\nu^r|+1})],\quad\xi\to 0.\ee
		If, in addition, $\phi$ has stable integer shifts, then we must have $\sm_\infty(a, 2I_d)>m$.
		
		\item[(2)] If $\sm_\infty(a,2I_d)>\ttm$, $\phi\in(C^m(\dR))^{r\times 1}$ and $\vgu_a$ satisfies \er{ghsd:vgu}, then the generalized Hermite subdivision scheme of type $\Lambda$ with the filter $a$ is convergent with limit functions in $C^m(\dR)$.
	\end{enumerate}	
\end{theorem}

\bp\begin{enumerate}
	
	\item[(1)] By letting $v=\td [\pe_j^r]^\tp$ with a limit function $\eta_v=\eta_{\td [\pe_j^r]^\tp}$ for $j\in\{1,\dots,r\}$ in \er{conv:ghsd}, we have
	\be\label{conv:ghsd:1}\lim_{n\to\infty}\|2^{|\nu^l|n}[\sd_{a, 2I_d}^n(\td I_r)]\pe_l^r-\partial^{\nu^l}\eta(2^{-n}\cdot)\|_{\dlrs{\infty}{r}{1}}=0,\quad\forall l=1,\dots,r,\ee
	where $\eta:=[\eta_{\td [\pe_1^r]^\tp},\dots,\eta_{\td [\pe_r^r]^\tp}]^\tp\in (C^m(\dR))^{r\times 1}$. We first show that $\wh{\eta}(0)\ne 0$. Since $\wh{\vgu_a}(0)\ne 0$, there must exist $k_0\in\{1,\dots,r\}$ such that $\wh{\vgu_a}(0)\pe_{k_0}^r\ne 0$. Also observe that \er{conv:ghsd:1} implies that \er{conv:vsd} holds with $\tau=\mu=\nu^{k_0}$. Hence, by item (1) of Theorem~\ref{thm:def:vsd}, we must have $\nu^{k_0}=0$ and $\wh{\vgu_a}(0)\wh{\eta}(0)=\wh{\vgu_a}(0)\pe_{k_0}^r\ne 0$, which, in particular, implies $\wh{\eta}(0)\ne 0$. 
	
	Next, we prove that $\eta=\wh{\vgu_a}(0)\pe_{k_0}^r\phi$, which implies $\phi\in (C^m(\dR))^{r\times 1}$. Using \er{conv:ghsd:1} with $l={k_0}$ and item (2) of Theorem~\ref{thm:def:vsd}, we deduce that $\wh{\eta}(2\xi)=\wh{a}(\xi)\wh{\eta}(\xi)$. By the condition \er{eigen:a:1}, there exists a unique $r\times 1$ vector of compactly supported distributions that satisfies \er{stan:vref} with $\dm=2$. Hence, as $\wh{\eta}(0)\ne 0$, we must have $\eta=c\phi$ for some $c\in\C\setminus\{0\}$. Moreover, the constant $c$ is determined by
	$$c=c\wh{\vgu_a}(0)\wh{\phi}(0)=\wh{\vgu_a}(0)\wh{\eta}(0)=\wh{\vgu_a}(0)\pe_{k_0}^r,$$
	and this proves the claim.

	Then, we prove \er{ghsd:vgu}. Since \er{conv:ghsd:1} holds, for every $l\in\{1,\dots,r\}$, we deduce from item (1) of Theorem~\ref{thm:def:vsd} that $\partial^\nu[\wh{\vgu_a}\pe^r_l](0)=0$ if $|\nu|<|\nu^l|$ and
	$$\partial^\nu[\wh{\vgu_a}\pe^r_l](0)=\partial^\nu[\wh{\vgu_a}\wh{\partial^{\nu^l}\eta}](0)=\td(\nu-\nu^l)i^{|\nu^l|}\nu^l!\wh{\vgu_a}(0)\wh{\eta}(0)\text{ if }|\nu|=|\nu^l|.$$
	Hence, for every $l\in\{1,\dots,r\}$, we have
	$$\wh{\vgu_a}(\xi)\pe^r_l=\wh{\vgu_a}(\xi)\wh{[\td \pe^r_l]}(\xi)=\beta_{\vgu_a,\td \pe^r_l,\nu^l}(i\xi)^{\nu^l}+\bo(\|\xi\|^{|\nu^l|+1}),\quad\xi\to 0,$$
	with
	$$\beta_{\vgu_a,\td \pe^r_l,\nu^l}=\frac{\partial^{\nu^l}[\wh{\vgu_a}\pe^r_l](0)}{i^{|\nu^l|}\nu^l!}=\wh{\vgu_a}(0)\wh{\eta}(0)=\wh{\vgu_a}(0)\pe^r_{k_0}.$$
	Therefore, up to multiplying a non-zero constant $\wh{\vgu_a}(0)\pe_{k_0}$, \er{ghsd:vgu} holds.

	Finally, we prove that $\sm_\infty(a,2I_d)>0$. By the definition of $\mom_{\vgu_a,\mu}$ in \er{mom}, it is clear that $\mom_{\vgu_a,0}=\dlrs{0}{r}{1}$. Then it is trivial that the set $G_{\vgu_a,0}:=\{\td \pe^r_1,\dots,\td \pe^r_r\}$ generates $\mom_{\vgu_a,0}$. For $l\in\{1,\dots,r\}$, if $\nu^l=0$, then using \er{ghsd:vgu} and \er{conv:ghsd:1}, we deduce that \er{conv:vsd:phi} must hold with $\dm =2$, $\mu=\nu^l=0$ and $u=\td \pe^r_l$. If $\nu^l\ne 0$, then \er{ghsd:vgu} yields $\beta_{\vgu_a,\td \pe_l^r,0}=0$. Using \er{conv:ghsd:1} and the fact that all entries of $\partial^{\nu^l}\eta$ are compactly supported and continuous, we deduce that
	\begin{align*}
	&\|[\sd_{a,2I_d}^n(\td I_r)]\pe^r_l-\beta_{\vgu_a,\td \pe^r_l,0}\phi(2^{-n}\cdot)\|_{\dlrs{\infty}{r}{1}}=\|[\sd_{a,2I_d}^n(\td I_r)]\pe^r_l\|_{\dlrs{\infty}{r}{1}}\\
	\le &2^{-|\nu^l|n}\|2^{|\nu^l|n}[\sd_{a,2I_d}^n(\td I_r)]\pe^r_l-\partial^{\nu^l}\eta(2^{-n}\cdot)\|_{\dlrs{\infty}{r}{1}}+2^{-|\nu^l|n}\|\partial^{\nu^l}\eta(2^{-n}\cdot)\|_{\dlrs{\infty}{r}{1}}\\
	\to &0,\quad n\to\infty.
	\end{align*} 
	Hence, \er{conv:vsd:phi} always holds with $\dm=2$, $\mu=0$ and $u=\td \pe^r_l$ ($l=1,\dots,r$). In other words, item (4) of Theorem~\ref{thm:conv:vsd} holds with $\dm =2$, $m=0$ and thus $\sm(a, 2I_d)>0$. If $\phi$ has stable integer shifts, then Theorem~\ref{thm:vcs} yields $\sm_\infty(a,2I_d)$.
	
	\item[(2)] By $\sm_\infty(a,2I_d)>\ttm$, \er{ghsd:vgu} and Theorem~\ref{thm:conv:vsd}, it is clear that \er{conv:ghsd:1} must hold. Now since, $\phi\in (C^m(\dR))^{r\times 1}$, using an argument similar to the proof of \er{conv:vsd:01} in item (5) of Theorem~\ref{thm:def:vsd}, we deduce from \er{conv:ghsd:1} that \er{conv:ghsd} must hold with $\eta_v:=v*\phi$ and this proves item (2).
\end{enumerate}

\ep

\begin{remark}Let us discuss how Theorem~\ref{thm:ghsd} improves the existing results \cite[Theorems 1.1 and 1.2]{han23-2} on generalized Hermite subdivision schemes.
	
	\begin{enumerate}
		\item[(1)]Suppose the generalized Hermite subdivision scheme with the filter $a$ is convergent with limit functions in $C^m(\dR)$. Under the additional assumption that $\phi$ satisfies
		\be\label{shift:phi}\spa\{\wh{\phi}(\pi\omega+2\pi k):\,k\in\dZ\}=\C^{r\times 1},\quad\forall \omega\in\Omega_{2I_d}=[0,1]^d\cap\dZ,\ee
		it has been shown in \cite[Theorem 1.1]{han23-2} that the filter $a$ must satisfy \er{eigen:a} with $\dm=2$, an order $m+1$ matching filter $\vgu_a$ is determined through \er{sr:diff:0} which, up to multiplying a non-zero constant, must be given by \er{ghsd:vgu}. Item (1) of Theorem~\ref{thm:ghsd} improves \cite[Theorem 1.1]{han23-2} by proving \er{sr:diff:0} under the weaker assumption \er{eigen:a:1} of the filter $a$.

		\item[(2)]Suppose $a$ has order $m+1$ sum rules with respect to $2I_d$ and $\vgu_a$ is an order $m+1$ matching filter that satisfies \er{ghsd:vgu}. Under the additional assumption that $\phi$ has stable integer shifts, \cite[Theorem 1.2]{han23-2} shows that $\sm_\infty(a,2I_d)>m$ must hold and it then follows from Theorem~\ref{thm:conv:vsd} that $\phi\in (C^m(\dR))^{r\times 1}$ and the generalized Hermite subdivision scheme with the filter $a$ is convergent with limit functions in $C^m(\dR)$. Without the stability on the integer shifts of $\phi$, item (2) proves the convergence of the generalized Hermit subdivision scheme with limit functions in $C^m(\dR)$ under the weaker assumptions $\sm_\infty(a,2 I_d)>\ttm$ and $\phi\in(C^m(\dR))^{r\times1}$. It is possible that $\phi\in (C^m(\dR))^{r\times 1}$ while $\sm_\infty(a,2I_d)<m$.

	\end{enumerate}	
\end{remark}

\section{Examples of Vector Subdivision Schemes}\label{sec:exmp}

Since the Lagrange and Hermite subdivision schemes have been extensively studied and many examples have been constructed in the literature, we will present vector subdivision schemes of other types in this section. In applications, the dilation matrix $2I_d$ is of the most interest in the literature of subdivision schemes and wavelet theory. In this section, we present some examples for the case $d=2$, that is, vector $2I_2$-subdivision schemes.

\subsection{Construction Guideline}

Let $\dm\in\N\setminus\{1\}$, $m\in\N_0$ and $r\in\N$. We first discuss how to construct a mask $a\in\dlrs{0}{r}{r}$ that yields a $C^m$-convergent vector $\dm I_d$-subdivision scheme. Here are the general construction steps:

\begin{enumerate}
	\item[(S1)] Parametrize the mask $a$ by
	$$\wh{a}(\xi)=\sum_{k\in [-N,N]^d\cap\dZ}a(k)e^{-ik\cdot\xi},\qquad\forall  \xi\in\dR,$$
	for some $N\in\N$ and some undetermined coefficients $a(k)$'s.
	
	\item[(S2)] \textbf{Sum rule conditions for $a$:} Parametrize a filter $\vgu_a\in\dlrs{0}{1}{r}$ by
	$$\wh{\vgu_a}(\xi)=\sum_{k\in [-K,K]^d\cap\dZ}\vgu_a(k)e^{-ik\cdot\xi},\qquad\forall  \xi\in\dR,$$
	for some $K\in\N$ and some undetermined coefficients $\vgu_a(k)$'s. Solve the linear system induced by \er{sr}, that is,
	$$\partial^\mu[\wh{\vgu_a}(\dm\cdot)\wh{a}](0)=\partial^\mu\wh{\vgu_a}(0),\quad \partial^\mu[\wh{\vgu_a}(\dm\cdot)\wh{a}](2\pi \omega)=0,\quad \forall \mu\in\cup_{q=0}\N_{0,q}^d,\quad \omega\in\Omega_{\dm I_d}\setminus\{0\},$$
	with the additional constraint $\wh{\vgu_a}(0)\ne 0$.
	Update $\vgu_a$ and $a$ by substituting the solutions of the above system.
	
	\item[(S3)]  \textbf{Try to optimize the $L_2$-smoothness exponent of $a$:} Choose the values of free parameters among the remaining free parameters such that $\sm_2(a,\dm I_d)$ is as large as possible. Ideally, try to achieve $\sm_2(a,\dm I_d)>m+\frac{d}{2}$ so that $\sm_\infty(a,\dm I_d)>m$. If not possible, try to directly estimate $\sm_\infty(a,\dm I_d)$ by using the structural properties of the mask $a$.
	
\end{enumerate}

If we add additional linear constraints to the above construction procedure, the mask $a$ can have a symmetry structure. The symmetry properties of multivariate filters are related to symmetry groups. By a \emph{symmetry group} in $\dR$, we mean a finite set $\mathcal{G}$ of $d\times d$ integer matrices that form a group under matrix multiplication. It is easy to see that any element $E$ of a symmetry group $\mathcal{G}$ must satisfy $\det(E)=\pm 1$. For any dimension $d$, there exists a \emph{trivial symmetry group} $\mathcal{G}=\{ I_d\}$, and any other symmetry groups are \emph{non-trivial}. Here are some examples of non-trivial symmetry groups:
\begin{itemize}
	
	\item For the univariate case $d=1$,  the only non-trivial symmetry group is $\mathcal{G}=\{-1,1\}$.
	
	\item For $d=2$, typical non-trivial symmetry groups include $\{-I_2,I_2\}$,
	\be\label{d4}D_4:=\left\{\pm I_2,\,\pm\begin{bmatrix}1 & 0\\
		0 & -1\end{bmatrix},\,\pm\begin{bmatrix}0 & 1\\
		1 & 0\end{bmatrix},\,\pm\begin{bmatrix}0 & 1\\
		-1 &0\end{bmatrix} \right\},\ee
	\be\label{d6}D_6:=\left\{\pm I_2,\,\pm\begin{bmatrix}0 & 1\\
		1 & 0\end{bmatrix},\,\pm\begin{bmatrix}-1 & 1\\
		0 & 1\end{bmatrix},\,\pm\begin{bmatrix}1 & 0\\
		1 & -1\end{bmatrix},\,\pm\begin{bmatrix}0 & 1\\
		-1 & 1\end{bmatrix},\, \pm\begin{bmatrix}1 & -1\\
		1 &0\end{bmatrix}\right\}.\ee
	$D_4$ is known as the \emph{full-axis symmetry group} and the symmetry group associated with the quadrilateral mesh in $\Z^2$. $D_6$ is known as the \emph{hexagon symmetry group} and the symmetry group associated with the triangular mesh in $\Z^2$.
	
\end{itemize}

Let $\phi=[\phi_1,\dots,\phi_r]^\tp$ be an $r\times 1$ vector of compactly supported distributions that satisfies \er{ref:eq:f} for some filter $a\in\dlrs{0}{r}{r}$. Let $T:=\{c_1,\dots,c_r\}$ be an ordered multiset with $c_1,\dots,c_r\in\dR$ and $\mathcal{G}$ be a symmetry group in $\dR$. We say that $\phi$ is \emph{$\cG$-symmetric about $T$ } if for every $E\in\mathcal{G}$, there exists an invertible $r\times r$ matrix $S_E=([S_E]_{l,t})_{1\le l,t\le r}$ such that
\be\label{sym:phi}\phi_l(E(x-c_l)+c_l)=\sum_{t=1}^r[S_E]_{l,t}\phi(x+c_t-c_l),\quad\forall x\in\dR,\quad l\in\{1,\dots,r\},\ee
or equivalently,
\be\label{sym:phi:f}\wh{\phi}(E^{-\tp}\xi)=D_T(E^{-\tp}\xi)S_ED_T(\xi)\phi(\xi),\quad\forall\xi\in\dR,\ee
where $D_T(\xi)$ is the $r\times r$ diagonal matrix given by
\be\label{sym:dt}D_T(\xi):=\DG(e^{ic_1\cdot\xi},\dots,e^{ic_r\cdot\xi}),\quad\forall\xi\in\dR.\ee
Using Theorem~\ref{thm:vcs} and the refinable equation \er{ref:eq:f}, one can derive that $\phi\in (C^m(\dR))^{r\times 1}$ is symmetric about $T$ with respect to $\mathcal{G}$ if and only if $\sm_\infty(a,\dm I_d)>m$ and $a$ is \emph{$\cG$-symmetric about $T$}, that is, 
\be\label{sym:a}\wh{a}(E^\tp\xi)=D_T(-\dm E^\tp\xi)S_{E^{-1}}D_T(\dm\xi)\wh{a}(\xi)D_T(-\xi)S_{E^{-1}}^{-1}D_T(E^\tp\xi),\quad\forall E\in\mathcal{G},\quad\xi\in\dR.\ee
If we require the $a$ in a vector $\dm I_d$-subdivision scheme to have symmetry, we first prescribed a multiset $T=\{c_1,\dots,c_r\}$ and a symmetry group $\mathcal{G}$, then add the linear constraint \er{sym:a} to the construction procedure that is given at the beginning of this subsection.

\subsection{Examples of Two-dimensional Scalar-type Vector Subdivision Schemes}

For $u\in\dlpp{0}$, define its \emph{filter support} by $\fsupp(u):=[k_1,k_2]\times[n_1,n_2]$ for some $k_1,k_2,n_1,n_2\in\Z$, which is the smallest two-dimensional interval such that $u(k)=0$ for all $k\in\fsupp(u)$. We use the following way to present a finitely supported filter $u\in l_0(\Z^2)$: suppose $\fsupp(u)=[k_1,k_2]\times[n_1,n_2]$, then we write
$$u=\begin{bmatrix}u(k_1,n_2) & u(k_1+1,n_2) &  \dots & u(k_2,n_2)\\
u(k_1,n_2-1) & u(k_1+1,n_2-1) &  \dots & u(k_2,n_2-1)\\
\vdots & \vdots &\ddots &\vdots\\
u(k_1, n_1) & u(k_1+1,n_1) &  \dots & u(k_2,n_1)\end{bmatrix}_{[k_1,k_2]\times[n_1,n_2]}.$$
For example, $\wh{u}(\xi_1,\xi_2)=e^{-i\xi_1}+2e^{i\xi_2}$ is presented as $u=\begin{bmatrix}
0 & 1\\
2 &0
\end{bmatrix}_{[0,1]\times[-1,0]}$.

\begin{exmp}\label{ex1}Let $\dm=d=r=2$, let $a=\begin{bmatrix}a_{11} & 0\\
a_{21} & a_{22}\end{bmatrix}\in(l_0(\Z^2))^{2\times 2}$ be the given by
$$a_{11}=\frac{1}{64}\begin{bmatrix}0 & 0 &  0 &  0 &  -1 &  -1 &  0\\
0 & 0 & -1 &  0 &  2 &  0 &  -1\\
0 &  -1 &  2 &  8 &  8 &  2 &  -1\\
0 &  0 & 8 &  16 &  8 &  0 &  0\\
-1 &  2 &  8 &  8 &  2 &  -1 &  0\\
-1 &  0 &  2 &  0 &  -1 &  0 &  0\\
0 &  -1 &  -1 &  0 &  0 &  0 &  0\end{bmatrix}_{[-3,3]^2},$$

$$a_{21}=\frac{1}{32}\begin{bmatrix}0 & 0 &  0 &  -1 &  0 &  0 &  -1\\
0 & 0 & 0 &  0 &  0 &  0 &  0\\
0 & 0 & 0 &  0 &  0 &  0 &  0\\
-1 & 0 & 0 &  0 &  0 &  0 &  -1\\
0 & 0 & 0 &  0 &  0 &  0 &  0\\
0 & 0 & 0 &  0 &  0 &  0 &  0\\
-1 &  0 &  0 &  -1 &  0 &  0 &  0\end{bmatrix}_{[-3,3]^2},\quad a_{22}=\frac{1}{64}\begin{bmatrix}0 & 0 & -1 & 0 & -1\\
0 & 0 & 0 & 0 & 0\\
-1 & 0 & 0 & 0 & -1\\
0 & 0 & 0 & 0 & 0\\
-1 & 0 & -1 & 0 & 0\end{bmatrix}_{[-2,2]^2}.$$
The filter $a$ is symmetric about $T=\{(0,0),\,(0,0)\}$ with respect to the group $D_6$. By calculation, we have $\sr(a, 2I_2)=4$ with an order $4$ matching filter given by $\wh{\vgu_a}(\xi)=[1,0]$. Moreover, $\sm_2(a,2I_2)\approx 2.4408$ and thus \er{sm:p2} yields $\sm_\infty(a, 2I_2)\ge 1.4408$. By Theorem~\ref{thm:conv:vsd}, the vector $2I_2$-subdivision scheme with the filter $a$ is at least $C^1$-convergent.

\end{exmp}

\begin{exmp}\label{ex2}Let $\dm=d=r=2$, let $a=\begin{bmatrix}a_{11} & 0\\
0 & a_{22}\end{bmatrix}\in(l_0(\Z^2))^{2\times 2}$ be the given by
$$a_{11}=\frac{1}{2048}\begin{bmatrix}0& 0 & 0 & 0 & 3 & 6 & 3 & 0 & 0 & 0 & 0\\
0 & 0 & 0 & 0 & 0 & 0 & 0 & 0 & 0 & 0 & 0\\
0 & 0 & 2 & 0 & -27 & -50 & -27 & 0 & 2 & 0 & 0\\
0 & 0 & 0 & 0 & 0 & 0 & 0 & 0 & 0 & 0 & 0\\
3 & 0 & -27 & 0 & 174 & 300 & 174 & 0 & -27 & 0 & 3\\
6 & 0 & -50 & 0 & 300 & 512 & 300 & 0 & -50 & 0 & 6\\
3 & 0 & -27 & 0 & 174 & 300 & 174 & 0 & -27 & 0 & 3\\
0 & 0 & 0 & 0 & 0 & 0 & 0 & 0 & 0 & 0 & 0\\
0 & 0 & 2 & 0 & -27 & -50 & -27 & 0 & 2 & 0 & 0\\
0 & 0 & 0 & 0 & 0 & 0 & 0 & 0 & 0 & 0 & 0\\
0& 0 & 0 & 0 & 3 & 6 & 3 & 0 & 0 & 0 & 0\end{bmatrix}_{[-5,5]^2},$$

$$a_{22}=\frac{1}{2048}\begin{bmatrix}0 & 0 &  -1 &  -2 &  -1 &  0 &  0\\
0 & 0 & 0 &  0 &  0 &  0 &  0\\
-1 &  0 &  0 &  0 &  0 &  0 &  -1\\
2 &  0 &  0 &  0 &  0 &  0 &  2\\
-1 &  0 &  0 &  0 &  0 &  0 &  -1\\
0 & 0 & 0 &  0 &  0 &  0 &  0\\
0 & 0 &  -1 &  -2 &  -1 &  0 &  0\end{bmatrix}_{[-3,3]^2}.$$
The filter $a$ is symmetric about $T=\{(0,0),\,(0,0)\}$ with respect to the group $D_4$. By calculation, we have $\sr(a, 2I_2)=6$ with an order $6$ matching filter given by $\wh{\vgu_a}(\xi)=[1,0]$. Moreover, $\sm_2(a,2I_2)\approx 3.1751$ and thus \er{sm:p2} yields $\sm_\infty(a, 2I_2)\ge 2.1751$. By Theorem~\ref{thm:conv:vsd}, the vector $2I_2$-subdivision scheme with the filter $a$ is at least $C^2$-convergent.

\end{exmp}

\subsection{Examples of Two-dimensional Balanced Vector Subdivision Schemes}

Let $\dm\in\N\setminus\{1\}$, $r\in\N$ and $m\in\N_0$. Let $N$ be a $d\times d$ integer matrix with $|\det(N)|=r$ and define
\be\label{ga:n}\Gamma_N:=\{\gamma_1,\dots,\gamma_r\}:=[N[0,1)^d]\cap\dZ\text{ with }\gamma_1=0.\ee
Suppose $a\in\dlrs{0}{r}{r}$ has sum rules of order $m+1$ with respect to $\dm I_d$ with a matching filter $\vgu_a\in\dlrs{0}{1}{r}$ that takes the form
\be\label{vgu:blc}\wh{\vgu_a}(\xi)=\wh{c}(\xi)[e^{iN^{-1}\gamma_1\cdot\xi},\dots,e^{iN^{-1}\gamma_r\cdot\xi}]+\bo(\|\xi\|^{m+1}),\quad\xi\to 0\ee
for some $c\in\dlp{0}$ with $\wh{c}(0)\ne 0$, then $a$ is said to be an \emph{order $m+1$ $E_N$-balanced filter}. Balanced filters play an important role in constructing multiwavelets and multiframelets with high balancing orders, which are of great interest in sparse discrete multiwavelet/multiframelet transforms. For detailed discussions on the theory of balanced filters and multiframelets, see \cite{han10,hl20pp,lu23} and many references therein.

A vector $\dm I_d$-subdivision scheme with an $E_N$-balanced filter is called an $E_N$-balanced vector $\dm I_d$-subdivision scheme. Using the dilation factor $2I_2$, we propose a general method to construct a two-dimensional balanced filter from a scalar refinement filter. We first recall the following result, a special case of \cite[Theorem 5.1]{han23-2}.

\begin{theorem}\label{thm:blc:flt}Let $\varphi$ be a $d$-variate compactly supported distribution such that $\wh{\varphi}(2\xi)=\wh{A}(\xi)\wh{\varphi}(\xi)$ for some $A\in\dlp{0}$ that satisfies $\wh{A}(0)=1$ and has order $m+1$ sum rules with respect to $2I_d$. Let $N$ be a $d\times d$ integer matrix such that $|\det(N)|=r\in\N$ and define $\Gamma_N:=\{\gamma_1,\dots,\gamma_r\}$ via \er{ga:n}.  Define
	$$\phi(x):=[\varphi(Nx-\gamma_1),\dots,\varphi(Nx-\gamma_r)]^\tp,\quad\forall x\in\dR,$$
	then $\phi$ is a $2I_d$-refinable vector function that satisfies $\wh{\phi}(2\xi)=\wh{a}(\xi)\wh{\phi}(\xi)$, where $a=(a_{jl})_{1\le j,l\le r}\in\dlrs{0}{r}{r}$ is given by
	\be\label{a:vec}a_{jl}(k)=A(Nk-2\gamma_j+\gamma_l),\quad\forall k\in\dZ,\quad j,l\in\{1,\dots,r\}.\ee 
	The filter $a$ satisfies $\sm_p(a,2I_d)=\sm_p(A,2I_d)$ for all $p\in[1,\infty]$. Moreover, $a$ is an order $m+1$ $E_N$-balanced filter with respect to $2I_d$ with a matching filter $\vgu_a\in\dlrs{0}{1}{r}$ that satisfies \er{vgu:blc} for some $c\in\dlp{0}$ such that $\wh{c}(0)\ne 0$ and $\wh{c}(2N^\tp\xi)\wh{A}(\xi)=\wh{c}(N^\tp\xi)+\bo(\|\xi\|^{m+1})$ as $\xi\to 0$. 	
\end{theorem}

\begin{rem}\label{rem:blc}In Theorem~\ref{thm:blc:flt}, if the filter $A$ has a symmetry structure and the matrix $N$ is properly chosen, then the filter $a$ defined via \er{a:vec} is likely to have symmetry properties as well. Let $\mathcal{G}$ be a non-trivial symmetry group of $d\times d$ integer matrices (non-trivial means $\mathcal{G}\ne\{I_d\}$). Suppose $A$ is $\mathcal{G}$-symmetric about $0$, that is,
	\be\label{sym:A}A(Ek)=A(k),\quad\forall k\in\dZ,\quad E\in\mathcal{G}.\ee
	Let $N$ be a $d\times d$ matrix such that $|\det(N)|=r$ and is compatible with $\cG$, that is,
	$$N^{-1}EN\in\cG,\quad \forall E\in\cG.$$
	Let $\gamma\in\dZ$ and define $a_\gamma(k):=A(Nk-\gamma)$ for all $k\in\dZ$. If $(I_d-E)N^{-1}\gamma\in\dZ$ for some $E\in\cG$, it is easy to deduce from \er{sym:A} that
	\be\label{sym:a:ga}a_\gamma(k)=a_\gamma(E(k-N^{-1}\gamma)+N^{-1}\gamma),\quad\forall k\in\dZ.\ee
	Therefore, if there exists a subgroup $\cH\subseteq\cG$ such that $(I_d-E)N^{-1}\gamma\in\dZ$ for all $E\in\cH$, then \er{sym:a:ga} holds for all $E\in\cH$, in other words, $a_{\gamma}$ is $\cH$-symmetric about $N^{-1}\gamma$. 
	
	Now let $\gamma=2\gamma_j-\gamma_l$ for $j,l\in\{1,\dots,r\}$, then $a_\gamma=a_{jk}$. Here we provide two examples to demonstrate the possible symmetry properties of $a_{jk}$:
	
	\begin{enumerate}
		
		\item[1.]Let $N:=M_{\sqrt{2}}:=\begin{bmatrix} 1 & 1\\
		1 &-1\end{bmatrix}$ be the quincunx dilation matrix. We have  $|\det(M_{\sqrt{2}})|=2$ and
		\be\label{ga:m:sq2}\Gamma_{M_{\sqrt{2}}}:=[M_{\sqrt{2}}[0,1)^2]\cap\Z^2:=\{\gamma_1:=(0,0),\, \gamma_2:=(1,0)\}.\ee
		By calculation, $M_{\sqrt{2}}$ is compatible with $\cG=D_4$ and $(I_2-E)M_{\sqrt{2}}^{-1}(2\gamma_j-\gamma_l)\in\Z^2$ for all $E\in D_4$ and $j,l\in\{1,2\}$.  Let $A\in l_0(\Z^2)$ be such that $\wh{A}(0,0)=1$ and is $D_4$-symmetric about $(0,0)$. Define $a=(a_{jl})_{1\le j,l\le 2}\in (l_0(\Z^2))^{2\times 2}$ by \er{a:vec} with $d=r=2$, then $a_{jl}$ is $D_4$-symmetric about $M_{\sqrt{2}}^{-1}(2\gamma_j-\gamma_l)$ for all $j,l\in\{1,2\}$, which, by further calculation, implies
		\be\label{sym:a:0}\wh{a}(E^\tp\xi)=D_T(-2E^\tp\xi)D(2\xi)\wh{a}(\xi)D_T(-\xi)D_T(E^\tp\xi),\quad\forall \xi\in\R^2,\quad E\in D_4,\ee
		where
		\be\label{t:q2}T:=\left\{M_{\sqrt{2}}^{-1}\gamma_1,\,M_{\sqrt{2}}^{-1}\gamma_2\right\}=\left\{(0,0),\,\left(\frac{1}{2},\frac{1}{2}\right)\right\},\ee
		$$D_T(\xi):=\DG(1, e^{\frac{i}{2}(\xi_1+\xi_2)}),\quad\forall\xi=(\xi_1,\xi_2)\in\R^2,$$
		that is, $a$ is $D_4$-symmetric about $T$.

		\item[2.]Let $N:=M_{\sqrt{3}}:=\begin{bmatrix} 1 & -2\\
		2 &-1\end{bmatrix}$ be the $\sqrt{3}$-dilation matrix. We have  $|\det(M_{\sqrt{3}})|=3$ and
		\be\label{ga:m:sq3}\Gamma_{M_{\sqrt{3}}}:=[M_{\sqrt{3}}[0,1)^2]\cap\Z^2:=\{\gamma_1:=(0,0),\, \gamma_2:=(-1,0),\,\gamma_3:=(0,1)\}.\ee
		By calculation, $M_{\sqrt{3}}$ is compatible with $\cG=D_6$ and $(I_2-E)M_{\sqrt{3}}^{-1}(2\gamma_j-\gamma_l)\in\Z^2$ for all $E\in \cH$ and $j,l\in\{1,2,3\}$, where $\cH\subseteq D_6$ is the following subgroup:
		\be\label{ch:d6}\cH:=\left\{I_2,\,\begin{bmatrix}0 &-1\\
			1 &-1\end{bmatrix},\,\begin{bmatrix}-1 & 1\\
			-1 &0\end{bmatrix},\,\begin{bmatrix}0 &-1\\
			-1 & 0\end{bmatrix},\,\begin{bmatrix}-1 & 1\\
			0 & 1\end{bmatrix},\,\begin{bmatrix}1 & 0\\
			1 &-1\end{bmatrix}\right\}.\ee
		Let $A\in l_0(\Z^2)$ be such that $\wh{A}(0,0)=1$ and is $D_6$-symmetric about $(0,0)$. Define $a=(a_{jl})_{1\le j,l\le 3}\in (l_0(\Z^2))^{3\times 3}$ by \er{a:vec} with $d=2$ and $r=3$, then $a_{jl}$ is $\cH$-symmetric about $M_{\sqrt{3}}^{-1}(2\gamma_j-\gamma_l)$ for all $j,l\in\{1,2,3\}$, which, by further calculation, implies
		\be\label{sym:a:1}\wh{a}(E^\tp\xi)=D_T(-2E^\tp\xi)D(2\xi)\wh{a}(\xi)D_T(-\xi)D_T(E^\tp\xi),\quad\forall \xi\in\R^2,\quad E\in \cH,\ee
		where
		\be\label{t:q3}T:=\left\{M_{\sqrt{3}}^{-1}\gamma_1,\,M_{\sqrt{3}}^{-1}\gamma_2,\,M_{\sqrt{3}}^{-1}\gamma_3\right\}=\left\{(0,0),\,\left(\frac{1}{3},\frac{2}{3}\right),\,\left(\frac{2}{3},\frac{1}{3}\right)\right\},\ee
		$$D_T(\xi):=\DG(1, e^{\frac{i}{3}(\xi_1+2\xi_2)},e^{\frac{i}{3}(2\xi_1+\xi_2)}),\quad\forall\xi=(\xi_1,\xi_2)\in\R^2,$$
		that is, $a$ is $\cH$-symmetric about $T$.
	\end{enumerate}

\end{rem}

With Theorem~\ref{thm:blc:flt} and Remark~\ref{rem:blc}, we construct two families of balanced vector $2I_2$-subdivision schemes with symmetry. 

\begin{exmp}\label{ex3}Let $\dm=2$ and $m\in\N_0$, the \emph{centred B-spline filter of order $2m$} is defined by
	$$\wh{a_{2m}^B}(\xi):=2^{-2m}(1+e^{-i\xi})^m(1+e^{i\xi})^m,\quad\forall \xi\in\dR.$$ 
	It is well-known that $\sr(a_{2m}^B,2)=2m$, $\sm_p(a_{2m}^B,2)=2m-1+\frac{1}{p}$ for all $p\in[1,\infty]$, and $a_{2m}^B$ satisfies the symmetry property $\wh{a_{2m}^B}(\xi)=\wh{a_{2m}^B}(-\xi)$. Now, define $A_{2m}^B\in l_0(\Z^2)$ via 
	\be\label{A:2m}\wh{A_{2m}^B}(\xi_1,\xi_2):=\wh{a_{2m}^B}(\xi_1)\wh{a_{2m}^B}(\xi_2),\quad \forall(\xi_1,\xi_2)\in\R^2.\ee
	Then we have $\sr(A_{2m}^B,2I_2)=\sr(a_{2m}^B,2)=2m$ and $\sm_p(A_{2m}^B,2I_2)=\sm_p(a_{2m}^B,2)=2m-1+\frac{1}{p}$ for all $p\in[1,\infty]$. Moreover, $A_{2m}^B$ is $D_4$-symmetric about $(0,0)$.
	
	Let $N:=M_{\sqrt{2}}:=\begin{bmatrix}1 & 1\\
	1 & -1\end{bmatrix}$ be the quincunx dilation matrix. Define $a_{2m}=([a_{2m}]_{jl})_{1\le j,l\le 2}\in(l_0(\Z^2))^{2\times 2}$ by
	$$[a_{2m}]_{jl}(k):=A_{2m}^B(M_{\sqrt{2}}k-2\gamma_j+\gamma_l),\quad\forall k\in\Z^2,\quad j,l\in\{1,2\},$$
	where $\Gamma_{M_{\sqrt{2}}}:=\{\gamma_1,\gamma_2\}$ is defined by \er{ga:m:sq2}. Then $\sr(a_{2m}, 2I_2)=\sr(A_{2m}^B,2I_2)=2m$ with an order $2m$ matching filter $\vgu_{2m}\in (l_0(\Z^2))^{1\times 2}$ given by
	$$\wh{\vgu_{2m}}=\wh{c_{2m}}(\xi)[e^{iM_{\sqrt{2}}^{-1}\gamma_1\cdot\xi},\,e^{iM_{\sqrt{2}}^{-1}\gamma_2\cdot\xi}]+\bo(\|\xi\|^{2m})=\wh{c_{2m}}(\xi)[1,\,e^{\frac{i}{2}(\xi_1+\xi_2)}]+\bo(\|\xi\|^{2m}),\quad\xi\to (0,0),$$
	where $c_{2m}\in l_0(\Z^2)$ satisfies $\wh{c_{2m}}(0)\ne 0$ and
	$$\wh{c_{2m}}(2M_{\sqrt{2}}^{\tp}\xi)\wh{A_{2m}^B}(\xi)=\wh{c_{2m}}(M_{\sqrt{2}}^{\tp}\xi)+\bo(\|\xi\|^{2m}),\quad\xi\to(0,0).$$
	Hence, $a_{2m}$ is an order $2m$ $E_{M_{\sqrt{2}}}$-balanced filter. By Theorem~\ref{thm:blc:flt}, we have $\sm_p(a_{2m},2I_2)=\sm_p(A_{2m}^B,2I_2)=2m-1+\frac{1}{p}$ for all $p\in[1,\infty]$. Furthermore, by Remark~\ref{rem:blc}, $a_{2m}$ is $D_4$-symmetric about $T$, where $T$ is given in \er{t:q2}.
	
	Here, we present two examples of $a_{2m}$:
	
	\begin{enumerate}
		\item[(1)] Let $m=2$, the filter $a_{4}\in (l_0(\Z^2))^{2\times 2}$ is given by
		$$[a_4]_{1,1}=\frac{1}{256}\begin{bmatrix}0& 0& 1& 0& 0\\
		0& 6& 16& 6& 0\\
		1& 16& 36& 16& 1\\
		0& 6& 16& 6& 0\\
		0& 0& 1& 0& 0\end{bmatrix}_{[-2,2]^2},\quad [a_4]_{1,2}=\frac{1}{64}\begin{bmatrix}0& 1& 1& 0\\
		1& 6& 6& 1\\
		1& 6& 6& 1\\
		0& 1& 1& 0\end{bmatrix}_{[-2,1]^2},$$
		$$[a_4]_{2,1}=\frac{1}{256}\begin{bmatrix}0& 0& 1& 0& 0\\
		0& 6& 16& 6& 0\\
		1& 16& 36& 16& 1\\
		0& 6& 16& 6& 0\\
		0& 0& 1& 0& 0\end{bmatrix}_{[-1,3]^2},\quad [a_4]_{2,2}=\frac{1}{64}\begin{bmatrix}0& 1& 1& 0\\
		1& 6& 6& 1\\
		1& 6& 6& 1\\
		0& 1& 1& 0\end{bmatrix}_{[-1,2]^2}.$$
		We have $\sr(a_4,2I_2)=4$ with a matching filter $\vgu_4$ given by
		$$\wh{\vgu_4}(\xi_1,\xi_2)=\wh{c_4}(\xi_1,\xi_2)[1,e^{\frac{i}{2}(\xi_1+\xi_2)}]+\bo(\|\xi\|^4),\quad(\xi_1,\xi_2)\to (0,0),$$
		where $c_4\in l_0(\Z^2)$ satisfies
		$$\wh{c_4}(\xi_1,\xi_2)=1+\frac{1}{12}(\xi_1^2+\xi_2^2)+\bo(\|\xi\|^4),\quad\xi\to (0,0).$$
		Moreover, $\sm_\infty(a_4,2I_2)=3$.Therefore, by Theorem~\ref{thm:conv:vsd}, the order $4$ $E_{M_{\sqrt{2}}}$-balanced vector $2I_2$-subdivision scheme with the filter $a_4$ is $C^2$-convergent.
		
		\item[(2)]Let $m=3$, the filter $a_{6}\in (l_0(\Z^2))^{2\times 2}$ is given by
		$$[a_6]_{1,1}=\frac{1}{4096}\begin{bmatrix}0& 0& 0& 1& 0& 0& 0\\
		0& 0& 15& 36& 15& 0& 0\\
		0& 15& 120& 225& 120& 15& 0\\
		1& 36& 225& 400& 225& 36& 1\\
		0& 15& 120& 225& 120& 15& 0\\
		0& 0& 15& 36& 15& 0& 0\\
		0& 0& 0& 1& 0& 0& 0\end{bmatrix}_{[-3,3]^2},\quad [a_6]_{1,2}=\frac{1}{2048}\begin{bmatrix}0&  0&  3&  3&  0&  0\\
		0&  10&  45&  45&  10&  0\\
		3&  45&  150&  150&  45&  3\\
		3&  45&  150&  150&  45&  3\\
		0&  10&  45&  45&  10&  0\\
		0&  0&  3&  3&  0&  0\end{bmatrix}_{[-3,2]^2},$$
		$$[a_6]_{2,1}=\frac{1}{4096}\begin{bmatrix}0& 0& 0& 1& 0& 0& 0\\
		0& 0& 15& 36& 15& 0& 0\\
		0& 15& 120& 225& 120& 15& 0\\
		1& 36& 225& 400& 225& 36& 1\\
		0& 15& 120& 225& 120& 15& 0\\
		0& 0& 15& 36& 15& 0& 0\\
		0& 0& 0& 1& 0& 0& 0\end{bmatrix}_{[-2,4]^2},\quad [a_6]_{2,2}=\frac{1}{2048}\begin{bmatrix}0& 0& 3& 3& 0& 0\\
		0& 10& 45& 45& 10& 0\\
		3& 45& 150& 150& 45& 3\\
		3& 45& 150& 150& 45& 3\\
		0& 10& 45& 45& 10& 0\\
		0& 0& 3& 3& 0& 0\end{bmatrix}_{[-2,3]^2}.$$
		We have $\sr(a_6,2I_2)=6$ with a matching filter $\vgu_6$ given by
		$$\wh{\vgu_6}(\xi_1,\xi_2)=\wh{c_6}(\xi_1,\xi_2)[1,e^{\frac{i}{2}(\xi_1+\xi_2)}]+\bo(\|\xi\|^6),\quad(\xi_1,\xi_2)\to (0,0),$$
		where $c_6\in l_0(\Z^2)$ satisfies
		$$\wh{c_6}(\xi_1,\xi_2)=1+\frac{1}{8}(\xi_1^2+\xi_2^2)+\frac{11}{640}\xi_1^2\xi_2^2+\frac{31}{3840}(\xi_1^4+\xi_2^4)+\bo(\|\xi\|^6),\quad\xi\to (0,0).$$
		Moreover, $\sm_\infty(a_6,2I_2)=5$. Therefore, by Theorem~\ref{thm:conv:vsd}, the order $6$ $E_{M_{\sqrt{2}}}$-balanced vector $2I_2$-subdivision scheme with the filter $a_6$ is $C^4$-convergent.
		
	\end{enumerate}

\end{exmp}

\begin{exmp}\label{ex4}Let $\dm=2$ and $m\in\N_0$, define the filter $u_m\in l_0(\Z^2)$ by
	\be\label{u:m}\wh{u_m}(\xi):=2^{-3m}(1+e^{-i\xi_1})^m (1+e^{-i\xi_2})^m (1+e^{i(\xi_1+\xi_2)})^m,\quad\forall \xi=(\xi_1,\xi_2)\in\R^2.\ee
	By calculation, $\sr(u_m,2)=2m$ and $u_m$ is $D_6$-symmetric about $(0,0)$. 
	
	Let $N:=M_{\sqrt{3}}:=\begin{bmatrix}1 & -2\\
	2 & -1\end{bmatrix}$ be the $\sqrt{3}$-dilation matrix. Define $a_{u,m}=([a_{u,m}]_{jl})_{1\le j,l\le 3}\in(l_0(\Z^2))^{3\times 3}$ by
	$$[a_{u,m}]_{jl}(k):=u_m(M_{\sqrt{3}}k-2\gamma_j+\gamma_l),\quad\forall k\in\Z^2,\quad j,l\in\{1,2,3\},$$
	where $\Gamma_{M_{\sqrt{3}}}:=\{\gamma_1,\gamma_2,\gamma_3\}$ is defined by \er{ga:m:sq3}. Then $\sr(a_{u,m}, 2I_2)=\sr(u_m,2I_2)=2m$ with an order $2m$ matching filter $\vgu_{u,m}\in (l_0(\Z^2))^{1\times 2}$ given by
	\begin{align*}\wh{\vgu_{u,m}}=&\wh{c_{u,m}}(\xi)[e^{iM_{\sqrt{3}}^{-1}\gamma_1\cdot\xi},\,e^{iM_{\sqrt{3}}^{-1}\gamma_2\cdot\xi},\,e^{iM_{\sqrt{3}}^{-1}\gamma_3\cdot\xi}]+\bo(\|\xi\|^{2m})\\
	=&\wh{c_{u,m}}(\xi)[1,\,e^{\frac{i}{3}(\xi_1+2\xi_2)},\,e^{\frac{i}{3}(2\xi_1+\xi_2)}]+\bo(\|\xi\|^{2m}),\quad\xi\to (0,0),\end{align*}
	where $c_{u,m}\in l_0(\Z^2)$ satisfies $\wh{c_{u,m}}(0)\ne 0$ and
	$$\wh{c_{u,m}}(2M_{\sqrt{3}}^{\tp}\xi)\wh{u_m}(\xi)=\wh{c_{u,m}}(M_{\sqrt{3}}^{\tp}\xi)+\bo(\|\xi\|^{2m}),\quad\xi\to(0,0).$$
	Hence, $a_{u,m}$ is an order $2m$ $E_{M_{\sqrt{3}}}$-balanced filter. By Theorem~\ref{thm:blc:flt}, we have $\sm_p(a_{u,m},2I_2)=\sm_p(u_m,2I_2)$ for all $p\in[1,\infty]$. Furthermore, by Remark~\ref{rem:blc}, $a_{u,m}$ is $\cH$-symmetric about $T$, where $\cH$ and $T$ are given by \er{ch:d6} and \er{t:q3} respectively.

	Here, we present two examples of $a_{u,m}$:
	
	\begin{enumerate}
		\item[(1)]Let $m=2$, the filter $a_{u,2}\in(l_0(\Z^2))^{3\times 3}$ is given by
		$$[a_{u,2}]_{11}=\frac{1}{32}\begin{bmatrix}0& 1& 1\\
		1& 5& 1\\
		1& 1& 0\end{bmatrix}_{[-1,1]^2},\, [a_{u,2}]_{12}=\frac{1}{64}\begin{bmatrix}1& 6& 1\\
		6& 6& 0\\
		1& 0& 0\end{bmatrix}_{[-1,1]\times[-2,0]},\, [a_{u,2}]_{13}=\frac{1}{64}\begin{bmatrix}0& 0& 1\\
		0& 6& 6\\
		1& 6& 1\end{bmatrix}_{[-2,0]\times[-1,1]},$$
		
		$$[a_{u,2}]_{21}=\frac{1}{64}\begin{bmatrix}1& 6& 1\\
		6& 6& 0\\
		1& 0& 0\end{bmatrix}_{[0,2]^2},\, [a_{u,2}]_{22}=\frac{1}{64}\begin{bmatrix}0& 0& 1\\
		0& 6& 6\\
		1& 6& 1\end{bmatrix}_{[-1,1]\times[0,2]},\, [a_{u,2}]_{23}=\frac{1}{32}\begin{bmatrix}0& 1& 1\\
		1& 5& 1\\
		1& 1& 0\end{bmatrix}_{[-1,1]\times[0,2]},$$	
		
		$$[a_{u,2}]_{31}=\frac{1}{64}\begin{bmatrix}0& 0& 1\\
		0& 6& 6\\
		1& 6& 1\end{bmatrix}_{[0,2]^2},\, [a_{u,2}]_{32}=\frac{1}{32}\begin{bmatrix}0& 1& 1\\
		1& 5& 1\\
		1& 1& 0\end{bmatrix}_{[0,2]\times[-1,1]},\, [a_{u,2}]_{33}=\frac{1}{64}\begin{bmatrix}1& 6& 1\\
		6& 6& 0\\
		1& 0& 0\end{bmatrix}_{[0,2]\times[-1,1]}.$$	
		We have $\sr(a_{u,2},2I_2)=4$ with a matching filter $\vgu_{u,2}$ given by
		$$\wh{\vgu_{u,2}}(\xi_1,\xi_2)=\wh{c_{u,2}}(\xi_1,\xi_2)[1,e^{\frac{i}{3}(\xi_1+2\xi_2)},e^{\frac{i}{3}(2\xi_1+\xi_2)}]+\bo(\|\xi\|^4),\quad(\xi_1,\xi_2)\to (0,0),$$
		where $c_{u,2}\in l_0(\Z^2)$ satisfies
		$$\wh{c_{u,2}}(\xi_1,\xi_2)=1+\frac{1}{6}(\xi_1^2+\xi_2\xi_2+\xi_2^2)+\bo(\|\xi\|^4),\quad\xi\to (0,0).$$
		Moreover, $\sm_2(a_{u,2},2I_2)=\sm_2(u_2,2I_2)\approx 3.5$ and thus $\sm_\infty(a_{u,2},2I_2)\ge 2.5$. Therefore, by Theorem~\ref{thm:conv:vsd}, the order $4$ $E_{M_{\sqrt{3}}}$-balanced vector $2I_2$-subdivision scheme with the filter $a_{u,2}$ is $C^2$-convergent.

		\item[(2)]Let $m=3$, the filter $a_{u,3}\in(l_0(\Z^2))^{3\times 3}$ is given by
		$$[a_{u,3}]_{11}=\frac{1}{512}\begin{bmatrix}0& 0& 0& 1& 0\\
		0& 1& 18& 18& 1\\
		0& 18& 56& 18& 0\\
		1& 18& 18& 1& 0\\
		0& 1& 0& 0& 0\end{bmatrix}_{[-2,2]^2},\quad [a_{u,3}]_{12}=\frac{3}{512}\begin{bmatrix}0& 0& 1& 1\\
		0& 4& 13& 4\\
		1& 13& 13& 1\\
		1& 4& 1& 0\end{bmatrix}_{[-2,1]^2},$$
		
		$$[a_{u,2}]_{13}=\frac{3}{512}\begin{bmatrix}0& 1& 4& 1\\
		1& 13& 13& 1\\
		4& 13& 4& 0\\
		1& 1& 0& 0\end{bmatrix}_{[-2,1]^2},\quad [a_{u,3}]_{21}=\frac{3}{512}\begin{bmatrix}0& 0& 1& 1\\
		0& 4& 13& 4\\
		1& 13& 13& 1\\
		1& 4& 1& 0\end{bmatrix}_{[-1,2]\times[0,3]},$$
		
		$$[a_{u,3}]_{22}=\frac{3}{512}\begin{bmatrix}0& 1& 4& 1\\
		1& 13& 13& 1\\
		4& 13& 4& 0\\
		1& 1& 0& 0\end{bmatrix}_{[-1,2]^2},\quad [a_{u,3}]_{23}=\frac{1}{512}\begin{bmatrix}0& 0& 0& 1& 0\\
		0& 1& 18& 18& 1\\
		0& 18& 56& 18& 0\\
		1& 18& 18& 1& 0\\
		0& 1& 0& 0& 0\end{bmatrix}_{[-2,2]\times[-1,3]},$$	
		
		$$[a_{u,3}]_{31}=\frac{3}{512}\begin{bmatrix}0& 1& 4& 1\\
		1& 13& 13& 1\\
		4& 13& 4& 0\\
		1& 1& 0& 0\end{bmatrix}_{[0,3]\times[-1,2]},\quad [a_{u,3}]_{32}=\frac{1}{512}\begin{bmatrix}0& 0& 0& 1& 0\\
		0& 1& 18& 18& 1\\
		0& 18& 56& 18& 0\\
		1& 18& 18& 1& 0\\
		0& 1& 0& 0& 0\end{bmatrix}_{[0,3]\times[-2,2]},$$
		
		$$[a_{u,3}]_{33}=\frac{3}{512}\begin{bmatrix}0& 0& 1& 1\\
		0& 4& 13& 4\\
		1& 13& 13& 1\\
		1& 4& 1& 0\end{bmatrix}_{[-1,2]^2}.$$	
		
		We have $\sr(a_{u,3},2I_2)=6$ with a matching filter $\vgu_{u,3}$ given by
		$$\wh{\vgu_{u,3}}(\xi_1,\xi_2)=\wh{c_{u,3}}(\xi_1,\xi_2)[1,e^{\frac{i}{3}(\xi_1+2\xi_2)},e^{\frac{i}{3}(2\xi_1+\xi_2)}]+\bo(\|\xi\|^6),\quad(\xi_1,\xi_2)\to (0,0),$$
		where $c_{u,3}\in l_0(\Z^2)$ satisfies
		$$\wh{c_{u,3}}(\xi_1,\xi_2)=1+\frac{1}{4}(\xi_1^2+\xi_2\xi_2+\xi_2^2)+\frac{1}{30}(\xi_1^4+2\xi_1^3\xi_2+3\xi_1^2\xi_2^2+2\xi_1\xi_2^3+\xi_2^4)+\bo(\|\xi\|^6),\quad\xi\to (0,0).$$
		Moreover, $\sm_2(a_{u,32},2I_2)=\sm_2(u_3,2I_2)\approx 5.5$ and thus $\sm_\infty(a_{u,3},2I_2)\ge 4.5$. Therefore, by Theorem~\ref{thm:conv:vsd}, the order $6$ $E_{M_{\sqrt{3}}}$-balanced vector $2I_2$-subdivision scheme with the filter $a_{u,3}$ is $C^4$-convergent.
		
	\end{enumerate}
	
\end{exmp}

\section*{Acknowledgements}
The author thanks Professor Bin Han from the University of Alberta for discussing several topics on subdivision schemes and sharing his C program routine on computing the $L_2$-smoothness exponents of matrix-valued filters.

\end{document}